\documentclass[report]{imsart}
\RequirePackage[OT1]{fontenc}
\RequirePackage{amsthm,amsmath}
\RequirePackage[numbers]{natbib}
\RequirePackage{hypernat}
\usepackage{amssymb,graphicx,url,bbm}
\usepackage[colorlinks,citecolor=blue,urlcolor=blue]{hyperref}

\startlocaldefs
\def\d{\delta}
\def\m{\mu}
\def\th{\theta}
\def\D{{\cal D}}
\def\P{{\mathbb P}}
\def\G{{\mathbb G}}
\def\H{{\mathbb H}}
\def\Q{{\mathbb Q}}

\def\a{\alpha}
\def\b{\beta}
\def\R{\mathbb R}
\def\dd{\Delta}

\def\d{\delta}

\def\D{{\cal D}}

\def\G{{\mathbb G}}
\def\H{{\mathbb H}}

\def\P{{\mathbb P}}

\def\labda1{\lambda_1}
\def\labda2{\lambda_2}
\def\m{\mu}

\def\e{\varepsilon}
\def\f{\phi}
\def\t{\tau}
\def\k{\kappa}

\def\s{\sigma}

\def\ss{\scriptscriptstyle}
\def\comment#1{\relax}

\def\=in{\mathop{\rm =}}

\newtheorem{theorem}{Theorem}[section]
\newtheorem{conjecture}{Conjecture}[section]

\newtheorem{lemma}{Lemma}[section]

\newtheorem{remark}{Remark}[section]

\numberwithin{equation}{section}
\theoremstyle{plain}
\endlocaldefs

\begin{document}

\begin{frontmatter}
\title{Nonparametric (smoothed) likelihood and integral equations}
\runtitle{M(S)LE and integral equations}

\begin{aug}
\author{\fnms{Piet} \snm{Groeneboom}\corref{}\ead[label=e1]{P.Groeneboom@tudelft.nl}
\ead[label=u1,url]{http://dutiosc.twi.tudelft.nl/\textasciitilde pietg/}}
\address{Delft Institute of Applied Mathematics,
Delft University of Technology,\\
Mekelweg 4, 2628 CD Delft,
The Netherlands\\
\printead{e1},
\printead{u1}
}
\affiliation{Delft University}
\runauthor{P.\ Groeneboom}
\end{aug}

\begin{abstract}
We show that there is an intimate connection between the theory of nonparametric (smoothed) maximum likelihood estimators for certain inverse problems and integral equations. This is illustrated by estimators for interval censoring and deconvolution problems. We also discuss the asymptotic efficiency of the MLE for smooth functionals in these models.
\end{abstract}



\end{frontmatter}


\section{Introduction}
\label{section:intro}
\setcounter{equation}{0}
There is an intimate connection between (nonparametric) maximum likelihood estimators for inverse problems and integral equations, a connection that does not seem to be well-known. In the present paper I will concentrate on maximum (smoothed) likelihood estimators for interval censored data, but maximum likelihood estimators for deconvolution will also be discussed. I will show that integral equations play a crucial role in the development of distribution theory for so-called ``smooth functionals" (of which moments are the prototype), based on the maximum likelihood estimator (MLE), but also in the development of the local limit theory of the MLE.

In \cite{piet_geurt_birgit:10} the maximum smoothed likelihood estimator (MSLE) was studied for the current status model, the simplest interval censoring model. It is called the interval censoring, case 1, model in \cite{Gr:91} and \cite{GrWe:92}. It was shown in \cite{piet_geurt_birgit:10} that, under certain regularity conditions,  the MSLE, evaluated at a fixed interior point,  converges at rate $n^{-2/5}$ to the real underlying distribution function, if one takes a bandwidth of order $n^{-1/5}$. This convergence rate is faster than the convergence rate of the non-smoothed maximum likelihood estimator, which is $n^{-1/3}$ in this situation, as shown in  \cite{Gr:91} and \cite{GrWe:92}. Moreover, the limit distribution is normal, in contrast with the limit distribution of the non-smoothed maximum likelihood estimator.

In the more realistic interval censoring model, there is an interval in which the relevant (unobservable) event takes place. This situation is in fact much more common, in particular in medical statistics. It is called the interval censoring, case 2, model in \cite{Gr:91} and \cite{GrWe:92}. In \cite{piet:12b} the local distribution theory for the MSLE was developed for this model and it was shown that, under a condition which is called the ``separation condition", the MSLE converges at rate $n^{-2/5}$ if the bandwidth is of the usual order $n^{-1/5}$ and that the MSLE has a normal limit distribution, again in contrast with the ordinary MLE. Here a (non-linear) integral equation plays again a crucial role.

It should be noted that one can also consider more observation times per hidden variable, for example the so-called ``case $k$" interval censoring model, where there are $k>2$ observation points per hidden variable, or the ``mixed case interval censoring model", considered in \cite{schick:00}, where there are a random number of observation times. From a computational point of view, however, the jump from the current status model to the interval censoring, case 2, model is the really big jump, since for current status data one can compute the MLE in one step, whereas an iterative algorithm is needed for the computation of the MLE in the interval censoring, case 2, model. Also, if there are $k>2$ observation points, only the interval containing the hidden variable will be relevant for the computation of the MLE, the other intervals can be discarded. In this sense the case 2 model is the most fundamental model, since it is the prototype for the ``case $k$" and ``mixed case" models, certainly from a computational point of view. Also, it is conjectured that the MLE for the case 2 model will attain a faster rate than the MLE in the current status model, see \cite{lucien:99} and \cite{piet_tom:11}, but such an increase in rate is not expected in going from case 2 to case $k>2$ or to mixed case interval censoring, unless one lets the number of observation times grow with the sample size  (see in particular \cite{lucien:99} for a discussion on the attainable minimax rates in these situations, which is not restricted to the use of maximum likelihood estimators). For these reasons I will concentrate on the interval censoring, case 2, model as a prototype for these models.

As noted in \cite{piet:12b}, the MSLE has the advantage over the MLE  and the smoothed maximum likelihood estimator (SMLE) that it can be used in situations where the MLE or SMLE cannot be used. For example, the MLE itself is proved to be inconsistent for the current status continuous mark model (see \cite{maathuis_wellner:08}), and the SMLE will inherit the bad properties of the MLE in this situation, and also be inconsistent. On the other hand, a version of the MSLE, based on histograms, is proved to be consistent for this model in \cite{piet_geurt_birgit:12b}. A similar phenomenon holds for the  two-dimensional right-censoring model, where the MLE is inconsistent (see \cite{tsai:86}) and the SMLE will not make this better. In this case the MSLE will, under appropriate smoothing of the observation distribution, also be consistent.

The MSLE can be viewed as an estimator minimizing a Kullback-Leibler distance and is therefore a natural generalization of the MLE, which minimizes the Kullback-Leibler distance of the distributions in the allowed class w.r.t.\ the unsmoothed empirical observation distribution (of course the Kullback-Leibler distance is not a real ``distance", but we follow the common convention of calling it a distance here). The difference between the MLE and the MSLE is that, in computing the MSLE, one starts by smoothing the empirical observation distribution of the data, and next looks for a distribution in the allowed class, closest to this smoothed observation distribution in Kullback-Leibler distance. In this way one can prevent the inconsistency properties of the MLE, as observed in \cite{maathuis_wellner:08} and \cite{tsai:86}, which have as common cause that the MLE tries to distribute mass on lower dimensional surfaces without using the surrounding information of the (higher dimensional) data.

We note here in passing the peculiar fact that in this Kullback-Leibler minimization, the part involving the minimization is ``on the left side" of the argument. In large deviation theory (in particular the large deviation theory associated with efficiency computations for test statistics), one usually has to deal with minimizing
$$
{\cal K}(Q,P)=\int \log\frac{dQ}{dP}\,dQ
$$
over $Q$ for fixed $P$, where ${\cal K}(Q,P)$ is the Kullback-Leibler distance between the probability measures $Q$ and $P$, see, e.g., \cite{piet:80} and \cite{GOR:79} (the Kullback-Leibler distance is infinite if $Q$ is not absolutely continuous with respect to $P$). But in the minimization needed to compute the M(S)LE, one has to minimize
$$
{\cal K}(Q,P)
$$
over $P$ for fixed $Q$, and this minimization problem is essentially different, and less theory is available. This has to do with the asymmetry of the Kullback-Leibler distance, which, for example, is not present with the Hellinger distance, which is a real distance.

We start in section \ref{section:smooth_functionalsIC} by studying smooth functionals for interval censoring, where we discuss the theory, developed in \cite{GeGr:96}, \cite{GeGr:97} and \cite{GeGr:99}. The notation, introduced here, will be used in the remainder of the paper. Section \ref{section:local_limit_IC} discusses a local limit result for interval censoring, the separated case, which was proved in \cite{piet:96} (see Theorem \ref{Th:local_limit_IC}), and discusses a conjecture for the non-separated case. Section \ref{section:local_limit_IC_MSLE} discusses a limit result recently proved in \cite{piet:12b} for the MSLE for interval censoring, showing that, under appropriate regularity conditions, the rate of the MLE can be improved to $n^{-2/5}$ by using the MSLE instead of the MLE. Also, the limit distribution is normal here, in contrast with the limit behavior of the MLE.

Section \ref{section:decon} takes a more heuristic turn, in the hope that researchers will pick up on this interesting topic, where there are still many open problems. It is based on ``Nachdiplom" lectures I gave at the ETH Z\"urich, in the fall of 2007, on the invitation of Sara van de Geer, and it is the first time these lectures appear (partly) in print. During these lectures, I tried to develop the theory of the integral equations, associated with deconvolution. The big hurdle here is the fact that the relevant efficient influence functions are unbounded near the edge of the domain on which they are defined. These functions can also only be numerically determined by solving the associated integral equations and do not have explicit representations, except in the case of uniform and exponential deconvolution. Nevertheless, pursuing this approach seems worthwhile, since there is little doubt in my mind that the MLE will automatically give efficient estimates of smooth functionals here, just as in the case of interval censoring, in contrast with the usual estimates, based on Fourier methods. Section \ref{section:decon_local} discusses results and conjectures for the local limit behavior of the MLE for deconvolution. Some of the conjectures go back more than 20 years, but have been proved for special cases in the mean time.

\section{Smooth functionals in the interval censoring model}
\label{section:smooth_functionalsIC}
\setcounter{equation}{0}
We recall the interval censoring, case 2, model. Let $X_1,\ldots,X_n$ be a sample of unobservable random variables
from an unknown distribution function $F_0$ on $[0,\infty)$.
Suppose that one can observe $n$ pairs $(T_i,U_i),$ independent of $X_i,$ where $U_i>T_i$. Moreover,
\begin{equation}
\label{def_indicators}
\Delta_{i1}\stackrel{\mbox{\small def}}=1_{\{X_i\le T_i\}},\,
\dd_{i2}\stackrel{\mbox{\small def}}=1_{\{T_i<X_i\le U_i\}}\mbox{ and }\dd_{i3}\stackrel{\mbox{\small def}}=1-\dd_{i1}-\dd_{i2},
\end{equation}
provide the only information one has on the position of the random variables
$X_i$ with respect to the observation times $T_i$ and $U_i$. In
this set-up one wants to estimate the unknown distribution function $F_0$, generating the ``unobservables" $X_i$.

If $F_0$ is an absolutely continuous distribution function, the MLE converges locally in distribution only at rate $n^{-1/3}$. But if one wants to estimate a so-called ``smooth functional", one can use the MLE to construct an asymptotically efficient estimate which converges at rate $n^{-1/2}$. Assuming that the density $f_0$, corresponding to the distribution function $F_0$, has support contained in an interval $[0,M]\subset[0,\infty)$, the smooth functionals of interest allow the following expansion:
$$
K(F)=K(F_0)+\int\k_{F_0}(x)\,d\left(F-F_0\right)(x)+O(\|F-F_0\|_2^2),\leqno\mbox{\rm (D1)}
$$
for  distribution functions $F$ with support contained in $[0,M]$, where $\|F-F_0\|_2$ is
the $L_2$-distance between the distribution function $F$ and $F_0$ w.r.t.~Lebesgue measure on
$[0,M]$. The function $\k_F$ is called the ``canonical gradient" of the functional $K$ w.r.t.\ $F$ (also called ``efficient influence function"), and is supposed to belong to the space $L_2^0(F)$, of square integrable functions $a$ w.r.t.\ the measure $dF$, satisfying $\int a\,dF=0$.  The derivative of the function $x\mapsto\kappa_F(x)$ (when it exists) will be denoted by $k_F(x)$.

Condition (D1) holds for a wider class of functionals than just the class of linear functionals.
For linear functionals
$$
K(F)=\int_0^M c(x)\,dF(x),
$$
we have ${\k}_F(x)=c(x)-\int c(x) dF(x)$, and (D1) even holds without the $\cal O$-term.
However, the functional 
$$
K(F)=\int F^2(x)\,w(x)\,dx
$$
where $w$ is a bounded weight function, has canonical gradient 
$$
\k_F(x)=2 \int_{s=x}^M F(s) \,w(s)\,ds-\int_{x=0}^M \,\left[\int_{s=x}^M 2\,F(s)
\,w(s)\,ds \right] \,dF(x)$$
and also satisfies (D1). For we have
\begin{align*}
&K(F)-K(F_0)=\int_0^M F(x)^2\,w(x)\,dx-\int_0^M F_0(x)^2\,w(x)\,dx\\
&=\int_0^M \left\{F(x)-F_0(x)\right\}^2\,w(x)\,dx+2\int_0^M F(x)F_0(x)w(x)\,dx-2\int_0^M F_0(x)^2\,w(x)\,dx\\
&=2\int_0^M \left\{F(x)-F_0(x)\right\}F_0(x)w(x)\,dx+O\left(\|F-F_0\|^2\right)\\
&=2\int_0^M \int_{s=x}^M F_0(s)w(s)\,ds\,d\left(F-F_0\right)(x)+O\left(\|F-F_0\|^2\right)\\
&=\int_0^M \k_{F_0}(x)\,d\left(F-F_0\right)(x)+O\left(\|F-F_0\|^2\right).
\end{align*}
The last equality holds since each constant integrates to zero w.r.t.\ $d(F-F_0)$ on the interval $[0,M]$. Here we use that the canonical gradient should belong to$L_2^0(F_0)$, an important property that sometimes seems to be overlooked in this kind of computation. In the next to last line of the display above we simply use integration by parts.

Condition (D1) suggests a recipe for proving efficiency and asymptotic normality of $K(\hat F_n)$, if $\hat F_n$ is the (ordinary unsmoothed) MLE. We first try to establish that
\begin{equation}
\label{L_2-bound}
\|\hat F_n-F_0\|^2=o_p\left(n^{-1/2}\right),
\end{equation}
and next try to prove, using the characterizing properties of the MLE,
\begin{equation}
\label{as_efficiency_relation1}
n^{1/2}\int\k_{F_0}(x)\,d\bigl(\hat F_n-F_0\bigr)(x)=n^{1/2}\int\th_{F_0}(t,u,\delta_1,\delta_2)\,d\left({\mathbb Q_n}-Q_0\right)(t,u,\delta_1,\delta_2)+o_p\left(n^{-1/2}\right),
\end{equation}
where $\th_{F_0}$ is given by:
\begin{equation}
\label{canon_gradientIC}
\th_{F_0}(t,u,\d_1,\d_2)=E\left\{\k_{F_0}(X)\bigm|(T_1,U_1,\dd_{11},\dd_{12})=(t,u,\d_1,\d_2)\right\}.
\end{equation}
It is at this point that the integral equations come into play, because $\th_{F_0}$ also has the representation
\begin{equation}
\label{thetaFrepresentation}
\th_{F_0}(t,u,\d_1,\d_2)=-\d_1\frac{\f_{F_0}(t)}{F_0(t)}-\d_2\,\frac{\f_{F_0}(u)-\f_{F_0}(t)}{F_0(u)-F_0(t)}
+\left(1-\d_1-\d_2\right)\frac{\f_{F_0}(u)}{1-F_0(u)}\,,
\end{equation}
where the function $\f_{F_0}$ is a solution of the integral equation (equation (6) on p.\ 77 of \cite{GeGr:96} and (5) on p.\ 204 of \cite{GeGr:97}):
\begin{equation}
\label{phi-equation}
\f_{F_0}(x)=d_{F_0}(x)\Bigl\{\frac{d}{dx}\left(\k_{F_0}(x)\right)-\int_{u=0}^x \dfrac{\f_{F_0}(x)-\f_{F_0}(u)}{F_0(x)-F_0(u)} \,g(u,x)\,du+
\int_{v=x}^{M} \dfrac{\f_{F_0}(v)-\f_{F_0}(x)}{F_0(v)-F_0(x)} \,g(x,v)\,dv\Bigr\},
\end{equation}
and where $d_{F_0}(x)$ is given by
$$
d_{F_0}(x)=\dfrac{F_0(x)\{1-F_0(x)\}}{g_1(x)\{1-F_0(x)\}+g_2(x)F_0(x)}\,.
$$
Here the indicators $\d_k$, $k=1,2$,  correspond to the indicators $\Delta_{ik}$ in (\ref{def_indicators}); $\mathbb Q_n$ is the empirical measure of the observations $(T_i,U_i,\dd_{i1},\dd_{i2})$, $i=1,\dots,n$, and the corresponding underlying probability measure is $Q_0$. The observation times $(T_i,U_i)$ have density $g$, with first marginal $g_1$ and second marginal $g_2$. How one gets from the canonical gradient $\k_{F_0}$ in the hidden space to the canonical gradient $\th_{F_0}$ in the observation space is further explained in \cite{GeGr:96}, \cite{GeGr:97} and \cite{piet:96}, where also the connection with theory, developed in \cite{vaart:91}, is explained. 

If one succeeds in proving relation (\ref{as_efficiency_relation1}), one has in one stroke established both asymptotic normality and asymptotic efficiency of $K(\hat F_n)$. In fact, one then gets, also using (\ref{L_2-bound}),
\begin{equation}
\label{as_efficiency_relation2}
n^{1/2}\left\{K(\hat F_n)-K(F_0)\right\}\stackrel{{\cal D}}\longrightarrow N(0,\s^2_{Q_0}),
\end{equation}
where $N(0,\s^2_{Q_0})$ is a univariate normal distribution with expectation zero and variance
$$
\s^2_{Q_0}=\int\th_{F_0}\bigl(t,u,\delta_1,\delta_2\bigr)^2\,dQ_0(t,u,\delta_1,\delta_2),
$$
which in terms of $\f_{F_0}$ becomes:
\begin{equation}
\label{limit_var}
\s^2_{Q_0}=\int\frac{\f_{F_0}(t)^2}{F_0(t)}\,g_1(t)\,dt+
\int\frac{\left\{\f_{F_0}(u)-\f_{F_0}(t)\right\}^2}{F_0(u)-F_0(t)}\,g(t,u)\,dt\,du+\int\frac{\f_{F_0}(u)^2}{1-F_0(u)}\,g_2(u)\,du.
\end{equation}
The asymptotic variance $\s^2_{Q_0}$ is the smallest asymptotic variance any regular estimator can attain.

The limit result (\ref{as_efficiency_relation2}) is proved in \cite{GeGr:97} under the strict separation hypothesis
$$
\P\left\{U_i-V_i>\e\right\}=1,
$$
for some $\e>0$ and some additional regularity conditions, and for the case that the joint density of $(U_i,V_i)$ is positive on the diagonal (in which case there can be arbitrarily small observation intervals $(U_i,V_i)$) in \cite{GeGr:99}. We will only give some background to the relations (\ref{L_2-bound}) and (\ref{as_efficiency_relation1}) for the separated case here. To discuss this in a simple setting, satisfying the conditions for the validity of (\ref{as_efficiency_relation2}), we take as $F_0$ the uniform distribution function on $[0,1]$ and as $g$ the uniform density on the upper triangle of the unit square with vertices $(0,\e)$, $(0,1)$ and $(1-\e,1)$, where $\e\in(0,1)$. Moreover, we take
$$
\k_F(x)=x-\int u\,dF(u),
$$
so the smooth functional we want to estimate is just the first moment of the distribution of $X$. This means that the integral equation (\ref{phi-equation}) boils down to:
\begin{equation}
\label{phi-equation2}
\f_{F_0}(x)=d_{F_0}(x)\left\{1-\int_{u=0}^x \dfrac{\f_{F_0}(x)-\f_{F_0}(u)}{x-u} \,g(u,x)\,du+
\int_{v=x}^{1} \dfrac{\f_{F_0}(v)-\f_{F_0}(x)}{v-x} \,g(x,v)\,dv\right\},
\end{equation}
where
\begin{align}
\label{bd_condition_phi_equation}
g(t,u)=\frac2{(1-\e)^2}1_{\{u-t>\e\}}\,,\qquad g_1(t)=\frac{2(1-t-\e)}{(1-\e)^2}1_{[0,1-\e]}(t)\,,\qquad g_2(v)=\frac{2(u-\e)}{(1-\e)^2}1_{[\e,1]}(u)\,.
\end{align}
Even in this simple setting, the integral equation (\ref{phi-equation}) does not have a simple solution and we have to develop general theory to show that a solution exists.

\begin{figure}[!ht]
\begin{center}
\includegraphics[scale=0.5]{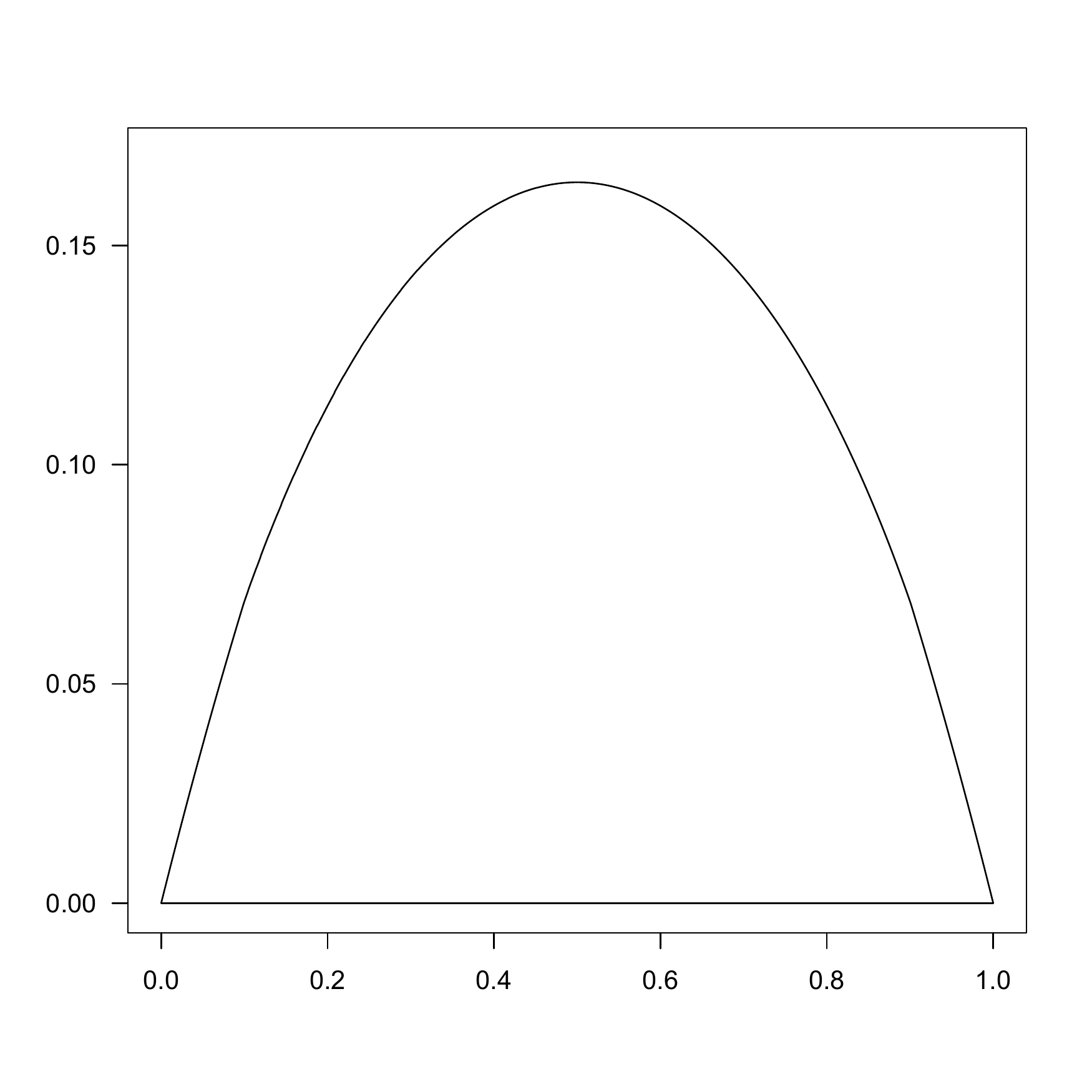}
\end{center}
\caption{The function $\f_{F_0}$, solving the integral equation (\ref{phi-equation2}), for $\e=0.1$.}
\label{fig:phi0}
\end{figure}

The integral equation (\ref{phi-equation2}) (and more generally, (\ref{phi-equation})), is a so-called Fredholm integral equation of the second kind, see, e.g. \cite{kress:89}. A picture of the solution $\f_{F_0}$ of (\ref{phi-equation2}) is shown in Figure \ref{fig:phi0}, where we took $\e=0.1$. By (\ref{thetaFrepresentation}) and (\ref{as_efficiency_relation2}), the asymptotic variance of
\begin{equation}
\label{sample_var}
n^{1/2}\int x\,d\bigl(\hat F_n-F_0\bigr)(x)
\end{equation}
is therefore given by (\ref{limit_var}), with $\f_{F_0}$ solving (\ref{phi-equation2}). Numerical solution of the integral equation (\ref{phi-equation2}) for $\e=0.1$ yielded
\begin{equation}
\label{sigma_Q_0}
\s^2_{Q_0}=\int\left\{\frac{\f_{F_0}(t)^2}{F_0(t)}+\frac{\left\{\f_{F_0}(u)-\f_{F_0}(t)\right\}^2}{F_0(u)-F_0(t)}+\frac{\f_{F_0}(u)^2}{1-F_0(u)}\right\}g(t,u)\,dt\,du
\approx0.11427,
\end{equation}
and a simulation study, using $10,000$ samples of size $n=1000$, yielded a variance $0.11470$ of the values of (\ref{sample_var}), so for sample size $n=1000$ the actual variance is close to the asymptotic variance, given by (\ref{limit_var}).

Because of the separation property $g(t,u)=0,\,u-t<\e$ and the hypothesis that $F_0$ has a continuous strictly positive density $f_0$ on its support (which, in the case of (\ref{phi-equation2}) is $[0,1]$), the integrating factors in the integrals on the right-hand side of (\ref{phi-equation}) and (\ref{phi-equation2}) are bounded. So we have a Fredholm integral equation of the second kind with a bounded integration kernel.

A further key to the treatment of the integral equation (\ref{phi-equation2}) is the following important observation. Assuming existence and uniqueness of the bounded continuous solution $\f_{F_0}$ of (\ref{phi-equation2}) (which is proved in \cite{GeGr:96}), we have:
\begin{equation}
\label{bound_phi}
\min_{x\in[0,1]}d_{F_0}(x)\le m\stackrel{\mbox{\small def}}=\min_{x\in[0,1]}\f_{F_0}(x)\le M\stackrel{\mbox{\small def}}=\max_{x\in[0,1]}\f_{F_0}(x)\le \max_{x\in[0,1]}d_{F_0}(x).
\end{equation}
The upper bound for the maximum $M$ follows in the following way. Let $x_0\in[0,1]$ be a point where the bounded continuous solution $\f_{F_0}$ attains its maximum $M$. Then
\begin{align*}
M=\f_{F_0}(x_0)&=d_{F_0}(x_0)\left\{1-\int_{u=0}^{x_0} \dfrac{M-\f_{F_0}(u)}{x_0-u} \,g(u,x_0)\,du-
\int_{v=x_0}^{1} \dfrac{M-\f_{F_0}(v)}{v-x_0} \,g(x_0,v)\,dv\right\}\\
&\le d_{F_0}(x_0)\le\max_{x\in[0,1]}d_{F_0}(x).
\end{align*}
The bound for $m$ is derived in a similar way. We similarly can obtain bounds for the derivative of $\f_{F_0}$ (using the bound we got for $\f_{F_0}$ itself).

Having established some properties of the solution $\f_{F_0}$ of the integral equation (\ref{phi-equation2}), we can explain why we can expect (\ref{L_2-bound}) and (\ref{as_efficiency_relation1}) to hold. First of all, we get
\begin{equation}
\label{L_2-bound1}
\|\hat F_n-F_0\|_2=O_p\left(n^{-1/3}\right).
\end{equation}
This can be proved by using some (by now) standard entropy methods, as developed, for example, in \cite{birman_solomjak:67} and \cite{ball_pajor:1990}.
Defining
$$
q_F(t,u,\d_1,\d_2)=\d_1 F(t)+\d_2\{F(u)-F(t)\}+(1-\d_1-\d_2)\{1-F(u)\},
$$
it is first proved, using \cite{birman_solomjak:67} (or \cite{ball_pajor:1990}), that the Hellinger distance $h(q_{\hat F_n},q_{F_0})$, defined by
$$
h(q_{\hat F_n},q_{F_0})=\left\{\tfrac12\int\left\{q_{\hat F_n}(t,u,\d_1,\d_2)^{1/2}-q_{F_0}(t,u,\d_1,\d_2)^{1/2}\right\}^2\,d\left({\mathbb Q}_n+Q_0\right)(t,u,\d_1,\d_2)\right\}^{1/2}
$$
satisfies
\begin{equation}
\label{hellinger_bound1}
h(q_{\hat F_n},q_{F_0})^2=O_p\left(n^{-2/3}\right),
\end{equation}
(part (i) of Corollary 2 in \cite{GeGr:97}), and next, using (\ref{hellinger_bound1}) and the inequalities
$$
\left(\hat F_n-F_0\right)^2\le 4\left(\sqrt{\hat F_n}-\sqrt{F_0}\right)^2
\mbox{ and }\left(\hat F_n-F_0\right)^2\le 4\left(\sqrt{1-\hat F_n}-\sqrt{1-F_0}\right)^2,
$$
that
$$
\bigl\|\hat F_n-F_0\bigr\|_2^2=O_p\left(n^{-2/3}\right).
$$
(part (ii) of Corollary 2 in \cite{GeGr:97}), which gives (\ref{L_2-bound1}).

Next we prove (\ref{as_efficiency_relation1}), using the following crucial lemma.

\begin{lemma}
\label{lemthF}
Let, in analogy with (\ref{thetaFrepresentation}), $\th_{\hat F_n}$ be defined by
\begin{equation}
\label{thetaFrepresentation_MLE}
\th_{\hat F_n}(t,u,\d_1,\d_2)=-\d_1\frac{\f_{\hat F_n}(t)}{\hat F_n(t)}-\d_2\,\frac{\f_{\hat F_n}(u)-\f_{\hat F_n}(t)}{\hat F_n(u)-\hat F_n(t)}
+\left(1-\d_1-\d_2\right)\frac{\f_{\hat F_n}(u)}{1-\hat F_n(u)}\,,
\end{equation}
where the function $\f_{\hat F_n}$ solves the integral equation
\begin{equation}
\label{phi-equation_MLE}
\f_{\hat F_n}(x)=d_{\hat F_n}(x)\left\{1-\int_{u=0}^x \dfrac{\f_{\hat F_n}(x)-\f_{\hat F_n}(u)}{\hat F_n(x)-\hat F_n(u)} \,g(u,x)\,du+
\int_{v=x}^{1} \dfrac{\f_{\hat F_n}(v)-\f_{\hat F_n}(x)}{\hat F_n(v)-\hat F_n(x)} \,g(x,v)\,dv\right\},
\end{equation}
and where $d_{\hat F_n}(x)$ is given by
$$
d_{\hat F_n}(x)=\dfrac{\hat F_n(x)\{1-\hat F_n(x)\}}{g_1(x)\{1-\hat F_n(x)\}+g_2(x)\hat F_n(x)}\,.
$$
Then
\begin{equation}
\label{obs_space_representation}
\int x\,d\bigl(\hat F_n-F_0\bigr)(x)=-\int \theta_{\hat F_n}(t,u,\d_1,\d_2)\,dQ_0(t,u,\d_1,\d_2).
\end{equation}
\end{lemma}

Lemma \ref{lemthF} is a special case of Lemma 1 on p.\ 214 of \cite{GeGr:97} and is the first step in proving (\ref{as_efficiency_relation1}). It gives a representation of the statistic of interest in the hidden space (the expression on the left side of (\ref{obs_space_representation})) in terms of a statistic in the observation space. The more general lemma in \cite{GeGr:97} holds for sufficiently well-behaved distribution functions $F$ instead of just $\hat F_n$ and the fact that $\hat F_n$ is the MLE is not used in the proof. We need, however, that 
 (\ref{phi-equation_MLE}) is a well-defined integral equation (at least for sufficiently large $n$, with probability tending to one), and this is not immediately clear. For example, the denominators of the integrands could be zero or arbitrarily close to zero. Moreover, $\hat F_n$ has jumps, so we have to deal with a mix of (absolutely) continuous functions and functions with jumps in this equation.

We now sketch the approach, taken in \cite{GeGr:97}. Let $J_i=[\t_i,\t_{i+1})$ be the intervals of constancy of $\hat F_n$, where $\t_0=0$ and $\t_{m+1}=1$. We define a piecewise constant version $\bar\f_{\hat F_n}$ of $\f_{\hat F_n}$ in the following way:
\begin{equation}
\label{def_bar_phi}
\bar\f_{\hat F_n}(x)=\left\{
\begin{array}{lll}
\f_{\hat F_n}(s),\,&\mbox{ if }\exists s\in J_i, \mbox{ such that }\hat F_n(s)=F_0(s),\\
\f_{\hat F_n}(\t_{i+1}-),\,&\mbox{ if }F_0(x)<\hat F_n(\t_i),\mbox{ for all }x\in J_i,\\
\f_{\hat F_n}(\t_i),\,&\mbox{ if }F_0(x)>\hat F_n(\t_i),\mbox{ for all }x\in J_i.
\end{array}
\right.
\end{equation}
We next define
\begin{equation}
\label{thetaFrepresentation_MLE3}
\bar\th_{\hat F_n}(t,u,\d_1,\d_2)=-\d_1\frac{\bar\f_{\hat F_n}(t)}{\hat F_n(t)}-\d_2\,\frac{\bar\f_{\hat F_n}(u)-\bar\f_{\hat F_n}(t)}{\hat F_n(u)-\hat F_n(t)}
+\left(1-\d_1-\d_2\right)\frac{\bar\f_{\hat F_n}(u)}{1-\hat F_n(u)}\,.
\end{equation}
Since $\bar\f_{\hat F_n}$ is absolutely continuous w.r.t.\ $\hat F_n$, we get:
$$
\int\bar\th_{\hat F_n}(t,u,\d_1,\d_2)\,d\Q_n(t,u,\d_1,\d_2)=0.
$$
Hence, using Lemma \ref{lemthF}, we get from (\ref{obs_space_representation}),
\begin{align*}
&\int x\,d\bigl(\hat F_n-F_0\bigr)(x)=-\int \theta_{\hat F_n}(t,u,\d_1,\d_2)\,dQ_0(t,u,\d_1,\d_2)\\
&=\int \bar\theta_{\hat F_n}(t,u,\d_1,\d_2)\,d\bigl(\Q_n-Q_0\bigr)(t,u,\d_1,\d_2)
+\int\left\{\bar\theta_{\hat F_n}(t,u,\d_1,\d_2)-\theta_{\hat F_n}(t,u,\d_1,\d_2)\right\}\,dQ_0(t,u,\d_1,\d_2).
\end{align*}
Empirical process theory yields:
$$
\sqrt{n}\int \bar\theta_{\hat F_n}(t,u,\d_1,\d_2)\,d\bigl(\Q_n-Q_0\bigr)(t,u,\d_1,\d_2)
\stackrel{{\cal D}}\longrightarrow N(0,\s^2_{Q_0}),
$$
where $\s^2_{Q_0}$ is defined by (\ref{sigma_Q_0}). Here we use
$$
\left\|\bar\f_{\hat F_n}-\f_{\hat F_n}\right\|_2+\left\|\f_{\hat F_n}-\f_{F_0}\right\|_2
+\left\|\hat F_n-F_0\right\|_2=O_p\left(n^{-1/3}\right),
$$
which is proved in \cite{GeGr:97}. In fact, defining
$$
\xi_{\hat F_n}=\frac{\f_{\hat F_n}}{\hat F_n(1-\hat F_n)}\,,\qquad\,\bar\xi_{\hat F_n}=\frac{\bar\f_{\hat F_n}}{\hat F_n(1-\hat F_n)}\,,
$$
it can be proved (see (31), p.\ 16 of \cite{GeGr:97} and Lemma 4 of \cite{GeGr:96}) that, for all $x\in[0,1]$,
$$
\left|\bar\xi_{\hat F_n}(x)-\xi_{\hat F_n}(x)\right|\le c_1\left|\hat F_n(x)-F_0(x)\right|\qquad\mbox{and}\qquad\left|\f_{\hat F_n}(x)-\f_{F_0}(x)\right|\le c_2\left|\hat F_n(x)-F_0(x)\right|,
$$
for positive constants $c_1$ and $c_2$. These properties of the function $\f_{\hat F_n}$ are derived from the integral equation (\ref{phi-equation_MLE}).

Note that $\f_{\hat F_n}$ has an absolutely continuous and a discrete part. Lemma 4 of \cite{GeGr:96} tells us that if $x$ and $y$ both belong to an interval $J_i$ between jumps, we have:
$$
\left|\f_{\hat F_n}(y)-\f_{\hat F_n}(x)\right|\le K_1|y-x|,
$$
for a positive constant $K_1$, independent of $J_i$, and that if $\hat F_n$ has a jump at $x$, we get:
$$
\left|\f_{\hat F_n}(x)-\f_{\hat F_n}(x-)\right|\le K_2|\hat F_n(x)-\hat F_n(x-)|
$$
for a positive constant $K_2$. So the discrete part of $\f_{\hat F_n}$ is absolutely continuous w.r.t.\ $\hat F_n$.

Finally, by Lemma 2, p.\ 215, of \cite{GeGr:97} we get:
$$
\int\left\{\bar\theta_{\hat F_n}(t,u,\d_1,\d_2)-\theta_{\hat F_n}(t,u,\d_1,\d_2)\right\}\,dQ_0(t,u,\d_1,\d_2)=O_p\left(n^{-2/3}\right).
$$
Going through steps of this type seems unavoidable if one wants to prove a result of type (\ref{as_efficiency_relation2}).

For the non-separated case a result of type (\ref{as_efficiency_relation2}) was proved in \cite{GeGr:99}, see Theorem 3.2 on p.\ 647 of \cite{GeGr:99}. The result is also discussed in \cite{piet:96}. In this case one can no longer use the integral equation (\ref{phi-equation_MLE}) directly because of the singularities of the integrand, but instead has to use a modified form of this integral equation in a transformed scale. The details are omitted here. A full discussion can be found in \cite{piet:96}.

\section{Local limit theory for the MLE in the interval censoring model}
\label{section:local_limit_IC}
\setcounter{equation}{0}
The distinction between the separated case ($\P\{U_i-V_i<\e\}=0$ for some $\e>0$) and the non-separated case, where we can have arbitrarily small observation intervals $(U_i,V_i)$, plays an even more prominent role in the local limit theory than in the theory for the smooth functionals.

For the non-separated case the following conjecture was launched in \cite{Gr:91} (and repeated in \cite{GrWe:92}).

\begin{conjecture}[Conjecture in \cite{Gr:91}]
\label{th:nonsep_conjecture}
Let $F_0$ and $H$ be continuously differentiable at $t_0$ and
$(t_0,t_0)$, respectively, with strictly positive derivatives
$f_0(t_0)$ and $h(t_0,t_0)$, where $H$ is the distribution function of
$(T_i,U_i)$. By continuous differentiability of $H$ at $(t_0,t_0)$ is
meant that the density $h(t,u)$ is continuous at $(t,u)$, if $t<u$ and
$(t,u)$ is sufficiently close to $(t_0,t_0)$, and that $h(t,t)$,
defined by
$$
h(t,t)=\lim_{u\downarrow t}h(t,u),
$$
is continuous at $t$, for $t$ in a neighborhood of $t_0$.

Let $0<F_0(t_0),H(t_0,t_0)<1$, and let $\hat F_n$ be the MLE of $F_0$. Then
$$
(n\log n)^{1/3}\left\{\hat F_n(t_0)-F_0(t_0)\right\}\bigm/\left\{\tfrac34f_0(t_0)^2/h(t_0,t_0)\right\}^{1/3}
\stackrel{{\cal D}}\longrightarrow 2Z,
$$
where $Z$ is the last time that standard two-sided Brownian motion minus the parabola $y(t)=t^2$ reaches its
maximum.
\end{conjecture}

It was also shown in \cite{Gr:91} that
Conjecture \ref{th:nonsep_conjecture} is true for a ``toy" estimator,
obtained by doing one step of the iterative convex minorant algorithm,
starting the iterations at the underlying distribution function $F_0$;
the ``toy" aspect is that we can of course not do this in practice. In
spite of the fact that now more than twenty years have passed since this
conjecture has been launched, it still has not been proved. 

For the separated case one can also introduce a toy estimator of the
same type and one can again formulate the ``working hypothesis" that
 the toy estimator and the MLE have the same pointwise limit
behavior. Anticipating that this would hold, the
asymptotic distribution of the toy estimator is derived in \cite{We:95} for the separated case,
under the following conditions.

\begin{enumerate}
\item[(C1)] The support of $F_0$ is an interval $[0,M]$, where $M<\infty$.
\item[(C2)] $F_0$ and $G$ have densities $f_0$ and $g$ w.r.t.~Lebesgue measure on $\R$ and
$\R^2$, respectively.
\item[(C3)] Let the functions $k_{1,\e}$ and $k_{2,\e}$ be defined by
$$
k_{1,\e}(u)=\int_u^M\dfrac{g(u,v)}{F_0(v)-F_0(u)}\,\{F_0(v)-F_0(u)<\e^{-1}\}\,dv,
$$
and
$$
k_{2,\e}(v)=\int_0^v\dfrac{g(u,v)}{F_0(v)-F_0(u)}\,\{F_0(v)-F_0(u)<\e^{-1}\}\,du.
$$
Then, for $i=1,2$ and each $\e>0$,
$$
\lim_{\a\to\infty}\a\int_{(t_0,t_0+t/\a]}k_i(u,\e\a)\,du=0.
$$
\item[(C4)] $0<F_0(t_0)<1$ and $0<H(t_0,t_0)<1$.
\end{enumerate}
The motivation for these conditions is given in \cite{We:95} and
actually becomes clear from the proof, which is not given here.

\begin{theorem}[\cite{We:95}]
\label{wellner}
Suppose that assumptions (C1) to (C4) hold. Let $k_i,\,i=1,2,$ be
defined by
$$
k_1(u)=\int_u^M\dfrac{g(u,v)}{F_0(v)-F_0(u)}\,dv,\mbox{ and }
k_2(v)=\int_0^v\dfrac{g(u,v)}{F_0(v)-F_0(u)}\,du,
$$
and suppose that $f_0,g_1,g_2,k_1$ and $k_2$ are continuous at $t_0$, where $g_1$ and $g_2$ are the
first and second marginal densities of $g$, respectively. Moreover, assume
$f_0(t_0)>0$. Then, if $F_n^{(1)}$ is the estimator of the distribution function $F_0$, obtained
after one step of the iterative convex minorant algorithm, starting the iterations with $F_0$,
we have
$$
n^{1/3}\{2\xi(t_0)/f_0(t_0)\}^{1/3}\{F_n^{(1)}(t_0)-F_0(t_0)\}\stackrel{\D}\longrightarrow 2Z,
$$
where $Z$ is the last time where standard two-sided Brownian motion minus the parabola
$y(t)=t^2$ reaches its maximum, and where
$$
\xi(t_0)=\frac{g_1(t_0)}{F_0(t_0)}+k_1(t_0)+k_2(t_0)+\frac{g_2(t_0)}{1-F_0(t_0)}.
$$
\end{theorem}

It is indeed proved in \cite{piet:96} that, under slightly stronger
conditions (the most important one being that an observation interval
always has length $>\e$, for some $\e>0$), the MLE has the same limit behavior, using the
same norming constants. The expression for the asymptotic variance in
the separated case is remarkably different from the conjectured
variance in the non-separated case, which only depends on $F_0$ via
$f_0(t_0)$, showing that only the local behavior, depending on the
density at $t_0$, is important for the asymptotic variance (assuming
that the working hypothesis holds).

Note that if $(T_i,U_i)$ is uniform on the upper triangle of the unit
square, with vertices $(0,\e)$, $(0,1)$ and $(1-\e,1)$, we have:
$$
g_1(u)=\frac{2(1-u-\e)}{(1-\e)^2}\,,\qquad g_2(v)=\frac{2(v-\e)}{(1-\e)^2}\,,
$$
see (\ref{bd_condition_phi_equation}), and, if $F_0$ is the uniform distribution function on $[0,1]$,
$$
k_1(u)=\frac{2\log\{(1-u)/\e\}}{(1-\e)^2}\,,\qquad k_2(v)=\frac{2\log(v/\e)}{(1-\e)^2}\,,
$$
so
$$
\xi(t_0)=\frac2{(1-\e)^2}\left\{\frac{1-t_0-\e}{t_0}+\log\left(\frac{t_0(1-t_0)}{\e^2}\right)+\frac{t_0-\e}{1-t_0}\right\}
$$
in this case.

Note that the scaling constants in the Conjecture \ref{th:nonsep_conjecture} and in Theorem \ref{wellner} are of a different order: in Conjecture \ref{th:nonsep_conjecture} the order is $(n\log n)^{-1/3}$ and in Theorem \ref{wellner} the order is $n^{-1/3}$ (which is also the order of convergence of the MLE for current status data). One of the reasons to believe that the rate of the MLE for the non-separated case is indeed of order $(n\log n)^{-1/3}$ is the fact that in \cite{lucien:99} a histogram-type estimator is constructed which locally achieves this rate. Moreover, a simulation study in \cite{piet_tom:11} which compares the MLE with the histogram estimator in \cite{lucien:99} shows that the MLE actually has a smaller variance than the histogram estimator for the cases analyzed there and has for large samples a variance which is close to the conjectured asymptotic variance. Nevertheless Conjecture \ref{th:nonsep_conjecture} still has to be proved, no doubt using the associated integral equations.

We shall now sketch how the integral equations enter into the proof of the local limit result for the separated case. As in the preceding section, we assume that $F_0$ is defined on $[0,M]$ and has a continuous derivative $f_0$, staying away from zero on $[0,M]$ (defining $f_0$ at the boundary points by its left and right limits). The proof starts by showing that
\begin{equation}
\label{uniform_sf_bound}
\sup_{t\in(0,M)}\sqrt{n}\int_0^t\bigl\{\hat F_n(u)-F_0(u)\}\,du=O_p(1),
\end{equation}
see Lemma 4.4 on p.\ 146 of \cite{piet:96}. This is done by studying the integral equation
\begin{align}
\label{st-flour_inteq}
\f(x)=d_F(x)\left\{1_{[0,t)}(x)-\int_0^x\frac{\f(x)-\f(u)}{F(x)-F(u)}g(u,x)\,du
+\int_x^M\frac{\f(u)-\f(x)}{F(u)-F(x)}g(x,u)\,du\right\},
\end{align}
using the same notation as in the preceding section, see, e.g., (\ref{phi-equation_MLE}). This equation has a right-continuous solution $\f_{t,F}$ for each $t\in[0,M]$, if we restrict the distribution functions $F$ to the set
$$
{\cal F}_{\d}=\left\{F\in{\cal F}_{[0,M]}^d:\sup_{x\in[0,M]}|F(x)-F_0(x)|\le\d\right\},
$$
where ${\cal F}_{[0,M]}^d$ is the set of discrete distribution functions on $[0,M]$ with finitely many points of jump, and where we choose $\d>0$ sufficiently small. Note that we may assume, with probability tending to one, that $\hat F_n$ belongs to ${\cal F}_{\d}$ for sufficiently large $n$.

According to Lemma 4.3 of \cite{piet:96}, the set of discontinuities of the solution $\f_{t,F}$ of the integral equation (\ref{st-flour_inteq}) is contained in the set of discontinuities of $F$, augmented by the point $t$ (which is the only jump of the function $1_{[0,t)}$ on $[0,M]$). Furthermore, again according to Lemma 4.3 of \cite{piet:96}, we get, if $x$ is a point of jump of $\f_{t,F}$,
$$
\left|\f_{t,F}(x)-\f_{t,F}(x-)\right|\le c\{F(x)-F(x-)+1\},
$$
for some $c>0$ only depending on $F_0$ and $\d$. For points $x<y$ in an interval not containing jumps of $F$ we have:
$$
\left|\f_{t,F}(y)-\f_{t,F}(x)\right|\le c'(y-x),
$$
for some constant $c'>0$, again only depending on $F_0$ and $\d$.

Defining, as before (see (\ref{thetaFrepresentation_MLE})):
$$
\th_{t,\hat F_n}(t,u,\d_1,\d_2)=-\d_1\frac{\f_{t,\hat F_n}(t)}{\hat F_n(t)}-\d_2\,\frac{\f_{t,\hat F_n}(u)-\f_{t,\hat F_n}(t)}{\hat F_n(u)-\hat F_n(t)}
+\left(1-\d_1-\d_2\right)\frac{\f_{t,\hat F_n}(u)}{1-\hat F_n(u)}\,,
$$
we get, following the proof of Lemma 4.4 in \cite{piet:96},
$$
\int_0^t\bigl\{\hat F_n(u)-F_0(u)\}\,du=\int\th_{t,\hat F_n}(t,u,\d_1,\d_2)\,dQ_0(t,u,\d_1,\d_2). 
$$
This reduces again the functional of interest to an integral in the observation space.
We would have the result (\ref{uniform_sf_bound}) if we could write:
$$
\int\th_{t,\hat F_n}(t,u,\d_1,\d_2)\,dQ_0(t,u,\d_1,\d_2)
=\int\th_{t,\hat F_n}(t,u,\d_1,\d_2)\,d\bigl(Q_0-\Q_n\bigr)(t,u,\d_1,\d_2).
$$
This would be true if $\f_{t,\hat F_n}$ would be absolutely continuous w.r.t.\ $\hat F_n$. But since this is not the case, we take a function $\bar\f_{t,\hat F_n}$ close to $\f_{t,\hat F_n}$ which is, apart from possibly having a jump at the point $t$, is absolutely continuous w.r.t.\ $\hat F_n$. If $t$ belongs to the interval $(\t_i,\t_{i+1})$ between successive jumps of $\hat F_n$,  $\bar\f_{t,\hat F_n}$ is defined by
$$
\bar\f_{t,\hat F_n}(x)=\left\{\begin{array}{lll}
\f_{t,\hat F_n}(t-),\,&\mbox{ if }x\in[\t_i,t),\\
\f_{t,\hat F_n}(t),\,&\mbox{ if }x\in[t,\t_{i+1}),
\end{array}
\right.
$$
on the other intervals $\bar\f_{t,\hat F_n}$ can be defined in the same way as in (\ref{def_bar_phi}).
Analogously to what we did before, we define $\bar\th_{t,\hat F_n}$ by:
$$
\bar\th_{t,\hat F_n}(t,u,\d_1,\d_2)=-\d_1\frac{\bar\f_{t,\hat F_n}(t)}{\hat F_n(t)}-\d_2\,\frac{\bar\f_{t,\hat F_n}(u)-\bar\f_{t,\hat F_n}(t)}{\hat F_n(u)-\hat F_n(t)}
+\left(1-\d_1-\d_2\right)\frac{\bar\f_{t,\hat F_n}(u)}{1-\hat F_n(u)}\,.
$$
Then:
$$
\left|\int\bigl\{\bar\th_{t,\hat F_n}-\th_{t,\hat F_n}\bigr\}\,dQ_0\right|\le\|\hat F_n-F_0\|_2^2=O_p\left(n^{-2/3}\right),
$$
(see p.\ 147 of \cite{piet:96}).

Returning to the functionals
$$
t\mapsto\psi_n(t)\stackrel{\mbox{\small def}}=\int_0^t\bigl\{\hat F_n(u)-F_0(u)\}\,du,
$$
we get that, if $x\mapsto \hat F_n(x)-F_0(x)$ is of constant sign on an interval $J_i=[\t_i,\t_{i+1})$,
$$
\sup_{u\in J_i}\left|\psi_n(u)\right|\le\max\left\{\left|\psi_n(\t_i)\right|,\left|\psi_n(\t_{i+1})\right|\right\},
$$
since the function $\psi_n$ is then either increasing or decreasing on $J_i$. If, on the other hand $F_0$ and $\hat F_n$ cross on the interval $J_i$, $\psi_n$ first increases and then decreases after the crossing point, noting that $\hat F_n$ is constant and that $F_0$ increases on $J_i$, so we get, if $t\in(\t_i,\t_{i+1})$,
\begin{align*}
\psi_n(\t_i)\wedge\psi_n(\t_{i+1})\le\psi_n(t)&=\int\th_{t,\hat F_n}\,dQ_0
=\int\bar\th_{t,\hat F_n}\,dQ_0+O_p\left(n^{-2/3}\right)\\
&=\int\bar\th_{t,\hat F_n}\,d\bigl(Q_0-\Q_n)+\int\bar\th_{t,\hat F_n}\,d\Q_n+O_p\left(n^{-2/3}\right)\\
&\le \int\bar\th_{t,\hat F_n}\,d\bigl(Q_0-\Q_n)+O_p\left(n^{-2/3}\right),
\end{align*}
where we use $\int\bar\th_{t,\hat F_n}\,d\Q_n\le0$,
which is a consequence of the so-called Fenchel duality conditions, characterizing the MLE (see (4.38) on p.\ 147 of \cite{piet:96}). So we have, apart from a remainder term of order $O_p(n^{-2/3})$, in all cases bounded $\psi_n(t)$ by the values of $\psi_n$ at the points $\t_i$ and an integral of the form
\begin{equation}
\label{emp_integral}
\int\bar\th_{t,\hat F_n}\,d\bigl(Q_0-\Q_n).
\end{equation}

But $\psi_n(\t_i)$ can be written
$$
\psi_n(\t_i)=\int\th_{\t_i,\hat F_n}\,dQ_0=\int\bar\th_{\t_i,\hat F_n}\,dQ_0+O_p\left(n^{-2/3}\right)
=\int\bar\th_{\t_i,\hat F_n}\,d\bigl(Q_0-\Q_n\bigr)+O_p\left(n^{-2/3}\right),
$$
using that for $t=\t_i$ the function $\f_{t,\hat F_n}$ is absolutely continuous w.r.t.\ $\hat F_n$ (the jump of the function $1_{[0,t)}$ is in that case at the same location as a jump of $\hat F_n$). So we have bounded $\psi_n(t)$ by the empirical integrals of the form (\ref{emp_integral}), and the result now follows by empirical process theory.

Having established (\ref{uniform_sf_bound}), we now also get for a class of functions $\cal G$ of right-continuous functions $g:[0,M]\to\R$ of uniformly bounded variation:
\begin{equation}
\label{unif_Donsker_bound}
\sup_{g\in{\cal G}}\left|\int_0^M g(x)\bigl\{\hat F_n(x)-F_0(x)\bigr\}\,dx\right|=O_p\left(n^{-1/2}\right),
\end{equation}
see Corollary 4.3 on p.\ 149 of \cite{piet:96}. Using (\ref{unif_Donsker_bound}) we can first of all establish:
$$
\sup_{x\in[0,M]}\big|\hat F_n(x)-F_0(x)\bigr|=O_p\left(n^{-1/4}\right),
$$
see Corollary 3.4 in \cite{piet:96}.

Next we observe that for an (open, closed or half-open) interval $J_n$, $\hat F_n$ satisfies
\begin{align}
\label{off-diagonal_expansion}
&\int_{t\in J_n}\left\{\frac{\d_2}{\hat F_n(u)-\hat F_n(t)}-\frac{\d_2}{F_0(u)-\hat F_n(t)}\right\}\,d\Q_n
=\int_{t\in J_n}\frac{\d_2\bigl\{\hat F_n(u)-F_0(u)\bigr\}}{\bigl\{\hat F_n(u)-\hat F_n(t)\bigr\}\bigl\{F_0(u)-\hat F_n(t)\bigr\}}\,d\Q_n\nonumber\\
&=\int_{t\in J_n}\frac{\{F_0(u)-F_0(t)\}\bigl\{\hat F_n(u)-F_0(u)\bigr\}}{\bigl\{\hat F_n(u)-\hat F_n(t)\bigr\}\bigl\{F_0(u)-\hat F_n(t)\bigr\}}\,g(t,u)\,du
+\int_{t\in J_n}\frac{\d_2\bigl\{\hat F_n(u)-F_0(u)\bigr\}}{\bigl\{\hat F_n(u)-\hat F_n(t)\bigr\}\bigl\{F_0(u)-\hat F_n(t)\bigr\}}\,d\bigl(\Q_n-Q_0\bigr).
\end{align}
For the first integral on the right-hand side of (\ref{off-diagonal_expansion}) we get from (\ref{unif_Donsker_bound}):
\begin{align*}
\int_{t\in J_n}\{F_0(u)-F_0(t)\}g_1(t)\left\{\int\frac{\bigl\{\hat F_n(u)-F_0(u)\bigr\}}{\bigl\{\hat F_n(u)-\hat F_n(t)\bigr\}\bigl\{F_0(u)-\hat F_n(t)\bigr\}}\,g(u|t)\,du\right\}\,dt
=O_p\left(n^{-1/2}\left|J_n\right|\right),
\end{align*}
where $|J_n|$ denotes the length of the interval $J_n$, and for the second integral on the right-hand side of (\ref{off-diagonal_expansion}) we have, if $|J_n|=O_p(n^{-1/4})$,
\begin{align*}
\int_{t\in J_n}\frac{\d_2\bigl\{\hat F_n(u)-F_0(u)\bigr\}}{\bigl\{\hat F_n(u)-\hat F_n(t)\bigr\}\bigl\{F_0(u)-\hat F_n(t)\bigr\}}\,d\bigl(\Q_n-Q_0\bigr)
=O_p\left(n^{-3/4}\right),
\end{align*}
implying that if $J_n$ is of order $O_p(n^{-1/4})$ both terms are of order $O_p(n^{-3/4})$.
Since
$$
\int_{u\in J_n}\left\{\frac{\d_2}{\hat F_n(u)-\hat F_n(t)}
-\frac{\d_2}{\hat F_n(u)-F_0(t)}\right\}\,d\Q_n
$$
can be treated in a similar way, we get:
\begin{align*}
&\int_{t\in J_n}\left\{\frac{\d_1}{\hat F_n(t)}-\frac{\d_2}{\hat F_n(u)-\hat F_n(t)}\right\}\,d\Q_n
+\int_{u\in J_n}\left\{\frac{\d_2}{{\hat F_n(u)-\hat F_n(t)}}
-\frac{1-\d_1-\d_2}{1-\hat F_n(u)}\right\}\,d\Q_n\\
&=\int_{t\in J_n}\left\{\frac{\d_1}{\hat F_n(t)}-\frac{\d_2}{F_0(u)-\hat F_n(t)}\right\}\,d\Q_n
+\int_{u\in J_n}\left\{\frac{\d_2}{{\hat F_n(u)-F_0(t)}}
-\frac{1-\d_1-\d_2}{1-\hat F_n(u)}\right\}\,d\Q_n+O_p\left(n^{-3/4}\right),
\end{align*}
if $\left|J_n\right|=O_p(n^{-1/4})$. Notice that this replaces the value of $\hat F_n$ in the ``off-diagonal" arguments of the integrand by the corresponding value of $F_0$.

Using this result, one can in fact derive the improved result
$$
\sup_{t\in[-c,c]}\left|\hat F_n(t_0+n^{-1/3}t)-F_0(t_0)\right|=O_p\left(n^{-1/3}\right),
$$
see Lemma 4.6 in \cite{piet:96}, for each $c>0$ and an interior point $t_0\in(0,M)$. This, in turn, means that if we let $J_n$ be an interval of order $O(n^{-1/3})$ around a fixed point $t_0\in(0,M)$, we get:
\begin{align*}
&\int_{t\in J_n}\left\{\frac{\d_1}{\hat F_n(t)}-\frac{\d_2}{\hat F_n(u)-\hat F_n(t)}\right\}\,d\Q_n
+\int_{u\in J_n}\left\{\frac{\d_2}{{\hat F_n(u)-\hat F_n(t)}}
-\frac{1-\d_1-\d_2}{1-\hat F_n(u)}\right\}\,d\Q_n\\
&=\int_{t\in J_n}\left\{\frac{\d_1}{F_0(t)}-\frac{\d_2}{F_0(u)-F_0(t)}\right\}\,d\Q_n
+\int_{u\in J_n}\left\{\frac{\d_2}{{F_0(u)-F_0(t)}}
-\frac{1-\d_1-\d_2}{1-F_0(u)}\right\}\,d\Q_n\\
&\qquad-\int_{t\in J_n}\bigl\{\hat F_n(t)-F_0(t)\bigr\}\left\{\frac{\d_1}{F_0(t)^2}+\frac{\d_2}{\{F_0(u)-F_0(t)\}^2}\right\}\,d\Q_n\\
&\qquad-\int_{u\in J_n}\bigl\{\hat F_n(u)-F_0(u)\bigr\}\left\{\frac{\d_2}{\{F_0(u)-F_0(t)\}^2}
+\frac{\d_3}{\{1-F_0(u)\}^2}\right\}\,d\Q_n+o_p\left(n^{-2/3}\right).
\end{align*}
In particular, if $J_n=[\t_n,v)$, where $\t_n$ is a point of jump of $\hat F_n$ such that $|\t_n-t_0|=O_p(n^{-1/3})$, and where $v=\t_n+n^{-1/3}w$, $w>0$, we get, by the characterization of the MLE:
\begin{align*}
&0\le n^{2/3}\int_{t\in J_n}\left\{\frac{\d_1}{\hat F_n(t)}-\frac{\d_2}{\hat F_n(u)-\hat F_n(t)}\right\}\,d\Q_n
+\int_{u\in J_n}\left\{\frac{\d_2}{{\hat F_n(u)-\hat F_n(t)}}
-\frac{1-\d_1-\d_2}{1-\hat F_n(u)}\right\}\,d\Q_n\\
&=n^{2/3}\left\{\int_{t\in J_n}\left\{\frac{\d_1}{F_0(t)}-\frac{\d_2}{F_0(u)-F_0(t)}\right\}\,d\Q_n
+\int_{u\in J_n}\left\{\frac{\d_2}{{F_0(u)-F_0(t)}}
-\frac{1-\d_1-\d_2}{1-F_0(u)}\right\}\,d\Q_n\right\}\\
&\qquad+n^{2/3}\int_{t\in J_n}\bigl\{F_0(t)-F_0(t_0)\bigr\}\left\{\frac1{F_0(t)}+\frac{1}{F_0(u)-F_0(t)}\right\}\,g(t,u)\,\,dt\,du\\
&\qquad+n^{2/3}\int_{u\in J_n}\bigl\{F_0(u)-F_0(t_0)\bigr\}\left\{\frac{1}{F_0(u)-F_0(t)}
+\frac{1}{1-F_0(u)}\right\}\,\,g(t,u)\,\,dt\,du\\
&\qquad-n^{2/3}\int_{t\in J_n}\bigl\{\hat F_n(t)-F_0(t_0)\bigr\}\left\{\frac{1}{F_0(t)}+\frac{1}{F_0(u)-F_0(t)}\right\}\,g(t,u)\,\,dt\,du\\
&\qquad-n^{2/3}\int_{u\in J_n}\bigl\{\hat F_n(u)-F_0(t_0)\bigr\}\left\{\frac{1}{F_0(u)-F_0(t)}
+\frac{1}{1-F_0(u)}\right\}\,g(t,u)\,\,dt\,du+o_p(1),
\end{align*}
where the {\it inequality} on the left becomes an {\it equality} if the right endpoint $v$ of $J_n$ is also a point of jump of $\hat F_n$. As a function of $w$ (in $\t_n+n^{-1/3}w$), the first term on the right-hand side converges to a Brownian motion process and the second and third term on the right-hand side converge to a parabolic drift added to this process. The last two terms converge to the greatest convex minorant of this Brownian motion plus parabolic drift process. So the MLE is indeed asymptotically equivalent to the toy estimator, given in Theorem \ref{wellner}, and its asymptotic distribution is therefore also given by Theorem \ref{wellner}. So we have the following result.

\begin{theorem}\mbox{\rm (Theorem 4.4 of \cite{piet:96}.)}
\label{Th:local_limit_IC}
Let the conditions 
\begin{enumerate}
\item[(i)] $g_1$ and $g_2$ are continuous, with $g_1(x)+g_2(x) > 0$ for all $x \in [0,M]$, 
\item[(ii)] $(u,v)\mapsto g(u,v)$ is continuous on its support, with uniformy bounded partial derivatives, except
at a finite number of points, where left and right (partial) derivatives exist,  
\item[(iii)] $\P\{V-U < \e_{\ss 0} \}=0 $ for some $\e_{\ss 0}$ with $0 < \e_{\ss 0} \le M/2$, so  $g$ does not have mass close to the diagonal,
\end{enumerate}
be satisfied and let $F_0$ be continuous with a bounded derivative
$f_0$ on $[0,M]$, satisfying
$$
f_0(x)\ge c>0,\,x\in(0,M),
$$
for some constant $c>0$. Then we have at each point $t_0\in(0,M)$:
$$
n^{1/3}\{2\xi(t_0)/f_0(t_0)\}^{1/3}\{\hat F_n(t_0)-F_0(t_0)\}\stackrel{\D}\longrightarrow 2Z,
$$
where $\xi$ and $Z$ are defined as in Theorem \ref{wellner}. Hence $\hat F_n$ has the same
asymptotic distribution as the toy estimator $F_n^{(1)}$ of Theorem \ref{wellner}.
\end{theorem}

\section{Local limit theory for the MSLE in the interval censoring model}
\label{section:local_limit_IC_MSLE}
\setcounter{equation}{0}
As mentioned in the introduction, the MLE $\hat F_n$ minimizes, as a function of $F$, the Kullback-Leibler distance
$$
{\cal K}(\Q_n,P_{n,F})=\int\frac{d\Q_n}{dP_{n,F}}\,d\Q_n,
$$
over probability measures $P_{n,F}$ in the allowed class, where $\Q_n$ is the empirical measure of the observations $(T_i,U_i,\dd_{i1},\dd_{i2})$, $i=1,\dots,n$, and $P_{n,F}$ is a measure, defined by
\begin{align}
\label{def_P_n}
&\int\psi(t,u,\d_1,\d_2)\,dP_{n,F}(t,u,\d_1,\d_2)\nonumber\\
&=\int\Bigl\{\psi(t,u,1,0)F(t)+\psi(t,u,0,1)\bigl\{F(u)-F(t)\bigr\}+\psi(t,u,0,0)\bigl\{1-F(u)\bigr\}\Bigr\}\,d\G_n(t,u),
\end{align}
for bounded measurable functions $\psi$ w.r.t.\ the product of the Borel $\s$-algebra on $\R_+^2$ and the counting measure on $\{(1,0),(0,1),(0,0)\}$,
and where $\G_n$ is the empirical distribution function of the observation pairs $(T_i,U_i)$.

On the other hand, the MSLE minimizes the Kullback-Leibler distance
\begin{equation}
\label{Kullback_MSLE}
{\cal K}(\tilde Q_n,\tilde P_{n,F})=\int \log\frac{d\tilde Q_n}{d\tilde P_{n,F}}\,d\tilde Q_n
\end{equation}
over $F$, where $\tilde Q_n$ is a smoothed version of $\Q_n$, defined by
\begin{align}
\label{def_tilde_Q_n}
&\int\psi(t,u,\d_1,\d_2)\,d\tilde Q_n(t,u,\d_1,\d_2)\nonumber\\
&=\int\psi(t,u,1,0)\,d\tilde Q_n(t,u,1,0)+\int\psi(t,u,0,1)\,d\tilde Q_n(t,u,1,0)+\int\psi(t,u,0,0)\,d\tilde Q_n(t,u,0,0),
\end{align}
and the three measures on the right-hand side are smoothed versions of the measures $\Q_n(t,u,1,0)$, $\Q_n(t,u,0,1)$ and $\Q_n(t,u,0,0)$, respectively. Furthermore, $\tilde P_{n,F}$ is defined by
\begin{align}
\label{def_tilde_P_n}
&\int\psi(t,u,\d_1,\d_2)\,d\tilde P_{n,F}(t,u,\d_1,\d_2)\nonumber\\
&=\int\Bigl\{\psi(t,u,1,0)F(t)+\psi(t,u,0,1)\bigl\{F(u)-F(t)\bigr\}+\psi(t,u,0,0)\bigl\{1-F(u)\bigr\}\Bigr\}\,d\tilde G_n(t,u),
\end{align}
where $d\tilde G_n$ is given by
$$
d\tilde G_n(t,u)=d\tilde Q_n(t,u,1,0)+d\tilde Q_n(t,u,0,1)+d\tilde Q_n(t,u,0,0).
$$
Minimizing (\ref{Kullback_MSLE}) is equivalent to maximizing the smoothed log likelihood
\begin{align}
\label{criterion_function1}
\ell(F)&=\int\log F(t)\,d\tilde Q_n(t,u,1,0)+\int \log\{F(u)-F(t)\}\,d\tilde Q_n(t,u,0,1)\nonumber\\
&\qquad\qquad\qquad\qquad\qquad\qquad\qquad\qquad+\int \log\{1-F(u)\}\,d\tilde Q_n(t,u,0,0)
\end{align}
over $F$, and the maximizing $F$, which we will denote by $\tilde F_n$, is called the MSLE.

We now give a more specific form of (\ref{criterion_function1}). Let, as before, $g$ be the joint density of the observation pairs $(T_i,U_i)$, with first marginal $g_1$ and second marginal $g_2$. Moreover, let the densities $h_{01}$, $h_{02}$ and $h_0$ be defined by
\begin{equation}
\label{def_h_0}
h_{01}(t)=F_0(t)g_1(t),\qquad h_{02}(u)=\{1-F_0(u)\}g_2(u),\qquad h_0(t,u)=\{F_0(u)-F_0(t)\}g(t,u).
\end{equation}
We define $\tilde h_{nj}$, $j=1,2$, and $\tilde h_n$ as the estimates of the densities $h_{0j}$, $j=1,2$, and the 2-dimensional density $h_0$, where
\begin{equation}
\label{def_h_nj}
\tilde h_{n1}(t)=\frac1n\sum_{i=1}^nK_{b_n}(t-T_i)\,\dd_{i1},\qquad \tilde h_{n2}(u)=\frac1n\sum_{i=1}^nK_{b_n}(u-U_i)\dd_{i3},
\end{equation}
\begin{equation}
\label{def_h_n}
\tilde h_n(t,u)=\frac1n\sum_{i=1}^nK_{b_n}(t-T_i)\,K_{b_n}(u-U_i)\dd_{i2},
\end{equation}
and
$$
K_{b_n}(x)=\frac1{b_n}K\left(\frac{x}{b_n}\right),
$$
for a symmetric continuously differentiable kernel $K$ with compact support, like the triweight kernel
\begin{equation}
\label{def_K}
K(x)=\tfrac{35}{32}\left(1-x^2\right)^31_{[-1,1]}(x),\,x\in\R.
\end{equation}
At points near the boundary we use a boundary correction by replacing the kernel $K$ by a linear combination of $K(u)$ and $uK(u)$. Details on the latter are given in \cite{piet:12b}. As in \cite{piet:12b}, we take $b_n\asymp n^{-1/5}$.

\begin{figure}[!ht]
\begin{center}
\includegraphics[scale=0.5]{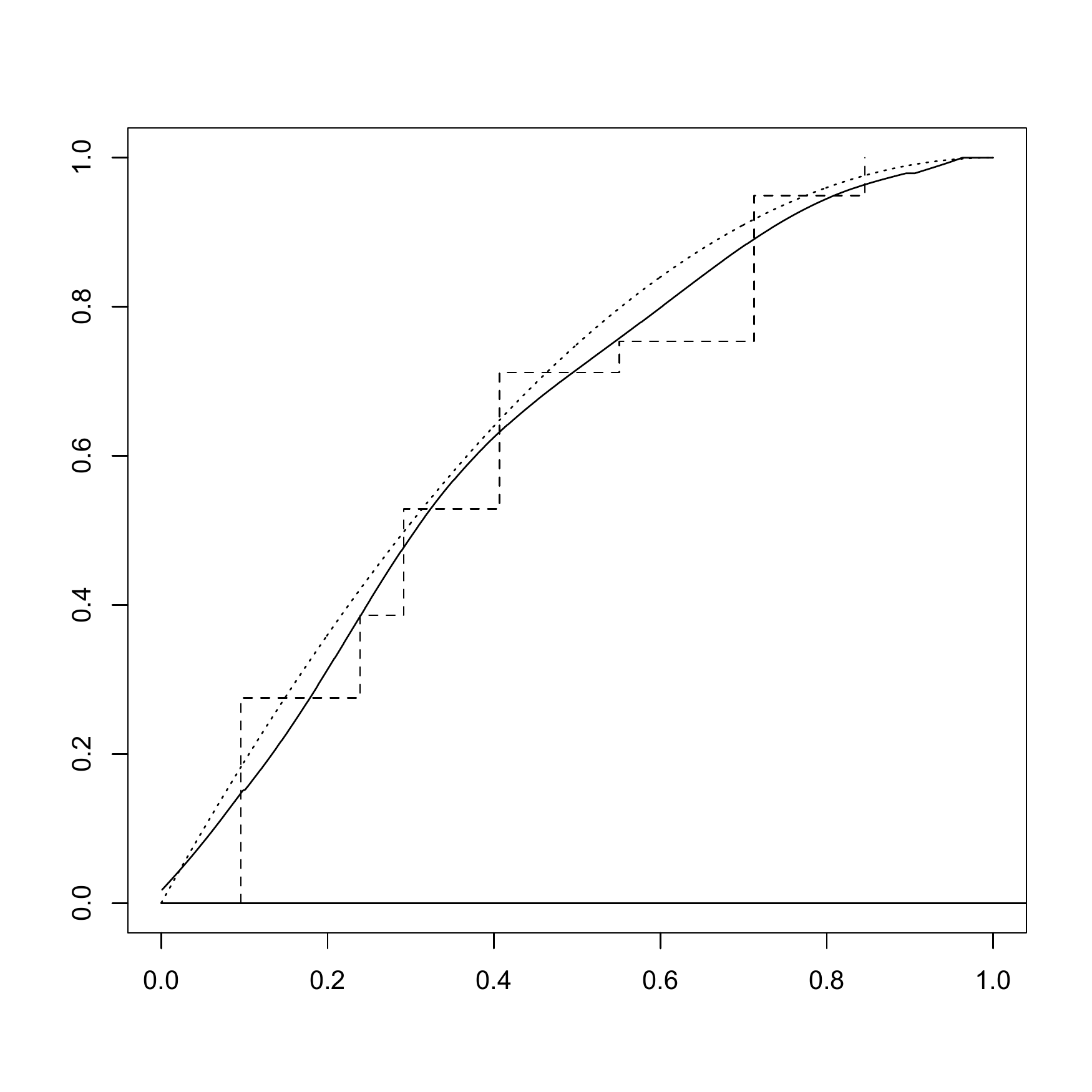}
\end{center}
\caption{The MSLE (solid) and MLE (dashed) on $[0,1]$ for a sample of size $n=100$ from the distribution function $F_0(x)=1-(1-x)^2$ (dotted); $g$ is uniform on the triangle with vertices $(0,\e)$, $(0,1)$ and $(1-\e,1)$, where $\e=0.1$. The bandwidth for the computation of the MSLE was $b_n=n^{-1/5}\approx0.398107$.}
\label{fig:MSLE100}
\end{figure}

\begin{figure}[!ht]
\begin{center}
\includegraphics[scale=0.5]{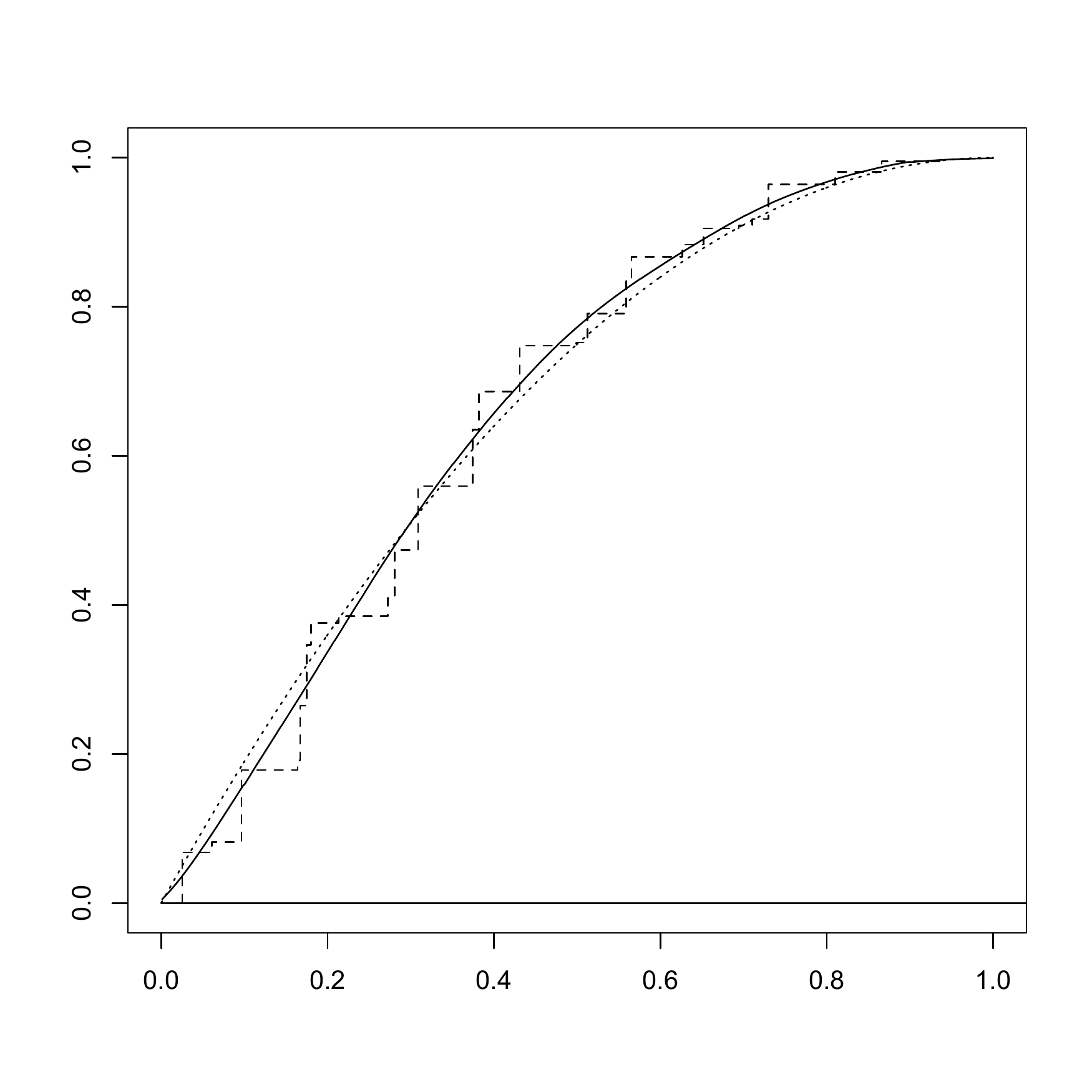}
\end{center}
\caption{The MSLE (solid) and MLE (dashed) on $[0,1]$ for a sample of size $n=1000$ from the distribution function $F_0(x)=1-(1-x)^2$ (dotted); $g$ is uniform on the triangle with vertices $(0,\e)$, $(0,1)$ and $(1-\e,1)$, where $\e=0.1$. The bandwidth for the computation of the MSLE was $b_n=n^{-1/5}\approx0.251189$.}
\label{fig:MSLE1000}
\end{figure}

With these definitions (\ref{criterion_function1}) takes the form:
\begin{align}
\label{criterion_function2}
\ell(F)=\int \tilde h_{n1}(t)\log F(t)\,dt+
\int \tilde h_{n2}(u)\{1-F(t)\}\,du+\int \tilde h_n(t,u)\log\{F(u)-F(t)\}\,dt\,du.
\end{align}
and the MSLE is the (sub-)distribution function, maximizing (\ref{criterion_function2}). Figures \ref{fig:MSLE100} and \ref{fig:MSLE1000} show the rather large improvement of the MSLE over the MLE when the underlying distribution is smooth. Note that, if we would not use boundary kernels, the measure $\tilde Q_n$ would be defined by
\begin{align*}
&\int\psi(t,u,\d_1,\d_2)\,d\tilde Q_n(t,u,\d_1,\d_2)\nonumber\\
&=\int\psi(t,u,1,0)\left\{\int K_{b_n}(t-x)K_{b_n}(u-y)\,d\Q_n(x,y,1,0)\right\}\,dt\,du\\
&\qquad+\int\psi(t,u,0,1)\left\{\int K_{b_n}(t-x)K_{b_n}(u-y)\,d\Q_n(x,y,0,1)\right\}\,dt\,du\\
&\qquad+\int\psi(t,u,0,0)\left\{\int K_{b_n}(t-x)K_{b_n}(u-y)\,d\Q_n(x,y,0,0)\right\}\,dt\,du.
\end{align*}
The analogous expression we obtain if boundary kernels are used, is obvious.

It is shown in \cite{piet:12b} that, under the separation hypothesis and some additional regularity conditions, the MSLE is asymptotically equivalent to the solution of a non-linear integral equation. We assume these conditions (given in Theorem 4.1 of \cite{piet:12b}) to be satisfied in the sequel. The relevant integral equation (in $F$) is given by:
\begin{align}
\label{fundamental_eq}
&\tilde h_{n1}(t)\{1-F(t)\}-\tilde h_{n2}(t)F(t)\nonumber\\
&\qquad\qquad+F(t)\{1-F(t)\}\left\{\int_{v=0}^t \frac{\tilde h_n(v,t)}{F(t)-F(v)}\,dv
-\int_{u=t}^M \frac{\tilde h_n(t,u)}{F(u)-F(t)}\,du\right\}=0,
\end{align}
see Lemma 4.5 of \cite{piet:12b}. Note that the corresponding equation for the underlying model:
\begin{align*}
&h_{01}(t)\{1-F(t)\}-h_{02}(t)F(t)\nonumber\\
&\qquad\qquad+F(t)\{1-F(t)\}\left\{\int_{v=0}^t \frac{h_0(v,t)}{F(t)-F(v)}\,dv
-\int_{u=t}^M \frac{h_0(t,u)}{F(u)-F(t)}\,du\right\}=0,
\end{align*}
is solved by $F_0$. Using the implicit function theorem in Banach spaces (\cite{dieudonne:1969}, Theorem 10.2.1), it is shown in \cite{piet:12b} that if $(\tilde h_{n1},\tilde h_{n2},\tilde h_n)$ is sufficiently close to $(h_{01},h_{02},h_0)$ in the supremum distance, equation (\ref{fundamental_eq}) has a unique solution $\tilde F_n$ in an open ball around $F_0$, again in the supremum distance. Next it is shown that $\tilde F_n$ coincides with the MSLE with probability tending to one and that
\begin{equation}
\label{first_bound_tilde_F_n}
\|\tilde F_n-F_0\|=O_p\left(n^{-3/10}\right),\,n\to\infty,
\end{equation}
where $\|\cdot\|$ denotes the supremum distance (part (ii) of Lemma 4.5 in \cite{piet:12b}). Note that starting with the non-sharp bound (\ref{first_bound_tilde_F_n}) is somewhat analogous to the approach in the derivation of the local limit behavior of the MLE in the preceding section, where first a bound on the supremum distance of order $O_p(n^{-1/4})$ was derived. In the derivation of (\ref{first_bound_tilde_F_n}) the implicit function theorem in Banach spaces is again used, but with a different norm for $(\tilde h_{n1},\tilde h_{n2},\tilde h_n)$ (instead of the supremum norm for $\tilde h_n$ a weaker integral-type norm is used).

Using the bound (\ref{first_bound_tilde_F_n}) it is subsequently shown that $\tilde F_n$ is close to the solution $\bar F_n$ of the {\it linear} integral equation
\begin{align}
\label{linear_asymp_equation}
&F(t)-F_0(t)
+d_{F_0}(t)\left\{\int_{u=0}^t\frac{g(u,t)\{F(t)-F_0(t)-F(u)+F_0(u)\}}{F_0(t)-F_0(u)}\,du\right.\nonumber\\
&\left.\qquad\qquad\qquad\qquad\qquad\qquad\qquad\qquad-\int_{u=t}^M\frac{g(t,u)\{F(u)-F_0(u)-F(t)+F_0(t)\}}{F_0(u)-F_0(t)}\,du\right\}\nonumber\\
&=\frac{\tilde h_{n1}(t)\{1-F_0(t)\}-\tilde h_{n2}(t)F_0(t)}{\{1-F_0(t)\}g_1(t)+F_0(t)g_2(t)}\nonumber\\
&\qquad\qquad\qquad\qquad+d_{F_0}(t)\left\{\int_{u=0}^t\frac{\tilde h_n(u,t)}{F_0(t)-F_0(u)}\,du-\int_{u=t}^M\frac{\tilde h_n(t,u)}{F_0(u)-F_0(t)}\,du\right\},
\end{align}
where $d_{F_0}$ is defined by
$$
d_{F_0}(t)=\frac{F_0(t)\{1-F_0(t)\}}{g_1(t)\{1-F_0(t)\}+g_2(t)F_0(t)}\,.
$$
In fact, it is shown that if $\bar F_n$ is the solution of (\ref{linear_asymp_equation}), we have:
$$
\|\tilde F_n-\bar F_n\|=O_p\left(n^{-3/5}\right),
$$
where $\|\cdot\|$ again denotes the supremum norm, which is a distance of smaller order than we can expect for the distance between $\tilde F_n$ and $F_0$.

The linear integral equation has properties which are analogous to the properties of the integral equations studied in the preceding sections, but is now an equation in $F$ itself instead of an equation in the associated function $\f_F$. In fact, an essential difference is that we now have asymptotic equalities and normality instead of asymptotic inequalities and non-normality. In the case of the MLE we had to infer the asymptotic properties via a functional of an associated process (greatest convex minorant of Brownian motion plus a parabolic drift), but we do not have to do this in the present case.

So $\bar F_n$ satisfies
\begin{align*}
&\bar F_n(t)-F_0(t)
+d_{F_0}(t)\left\{\int_{u=0}^t\frac{g(u,t)\{\bar F_n(t)-F_0(t)-\bar F_n(u)+F_0(u)\}}{F_0(t)-F_0(u)}\,du\right.\nonumber\\
&\left.\qquad\qquad\qquad\qquad\qquad\qquad\qquad\qquad-\int_{u=t}^M\frac{g(t,u)\{\bar F_n(u)-F_0(u)-\bar F_n(t)+F_0(t)\}}{F_0(u)-F_0(t)}\,du\right\}\nonumber\\
&=\frac{\tilde h_{n1}(t)\{1-F_0(t)\}-\tilde h_{n2}(t)F_0(t)}{\{1-F_0(t)\}g_1(t)+F_0(t)g_2(t)}\nonumber\\
&\qquad\qquad\qquad\qquad+d_{F_0}(t)\left\{\int_{u<t}\frac{\tilde h_n(u,t)}{F_0(t)-F_0(u)}\,du-\int_{u>v}\frac{\tilde h_n(t,u)}{F_0(u)-F_0(t)}\,du\right\}.
\end{align*}
As in the preceding section, the ``off-diagonal" terms $\bar F_n(u)-F_0(u)$ in the integrands on the left give a contribution of lower order $O_p(n^{-1/2})$ (``smooth functionals" again!) and we find that $\bar F_n$ is asymptotically equivalent to the ``toy estimator" $F_n^{toy}$, satisfying
\begin{align}
\label{toy-equation}
&\{F_n^{toy}(t)-F_0(t)\}\left\{1+d_{F_0}(t)\left\{\int_{u<t}\frac{g(u,t)}{F_0(t)-F_0(u)}\,dt
+\int_{u>v}\frac{g(t,u)}{F_0(u)-F_0(t)}\,du\right\}\right\}\nonumber\\
&=\frac{\tilde h_{n1}(t)\{1-F_0(t)\}-\tilde h_{n2}(t)F_0(t)}{\{1-F_0(t)\}g_1(t)+F_0(t)g_2(t)}\nonumber\\
&\qquad\qquad\qquad\qquad\qquad+d_{F_0}(t)\left\{\int_{u<t}\frac{\tilde h_n(u,t)}{F_0(t)-F_0(u)}\,du-\int_{u>v}\frac{\tilde h_n(t,u)}{F_0(u)-F_0(t)}\,du\right\}.
\end{align}
Since, by standard theory for kernel estimators, the right-hand side, multiplied by $n^{2/5}$, converges in distribution to a normal distribution, we get that also $n^{2/5}\{\hat F_n(t)-F_0(t)\}$ converges to a normal distribution (the same limit distribution as that of $n^{2/5}\{F_n^{toy}(t)-F_0(t)\}$), where $\hat F_n$ denotes the MSLE. The full result, with explicit asymptotic bias and  variance, is given below.

\begin{theorem}
\label{th:asymp_nor_sep}
Let Condition (1.1) of \cite{piet:12b} be satisfied. Moreover,
let $F_0$ be twice differentiable, with a bounded continuous derivative $f_0$ on the interior of $[0,M]$, which is bounded away from zero on $[0,M]$, with a finite positive right limit at $0$ and a positive left limit at $M$. Also, let $f_0$ have a bounded continuous derivative on $(0,M)$ and let $g_1$ and $g_2$ be twice differentiable on the interior of their supports $S_1$ and $S_2$, respectively, and let $g_1\{1-F_0\}+g_2F_0$ stay away from zero on $[0,M]$, where $g_1$ and $g_2$ are the marginals of the joint density $g$ of the pair of observation times $(T_i,U_i)$.

Furthermore, let $g$ have a bounded (total) second derivative on the interior of its support $S$, having finite limits approaching the boundary of $S$. Suppose that $X_i$ is independent of $(T_i,U_i)$, and let $d_{F_0}$ be defined by
$$
d_{F_0}(v)=\frac{F_0(v)\{1-F_0(v)\}}{g_1(v)\{1-F_0(v)\}+F_0(v)g_2(v)}\,.
$$
Then, if $b_n\asymp n^{-1/5}$, we have for each $v\in(0,M)$, and the MSLE $\hat F_n$,
$$
\sqrt{nb_n}\left\{\hat F_n(v)-F_0(v)-\frac{\b(v)b_n^2}{2\s_1(v)}\right\}\stackrel{{\cal D}}\longrightarrow
N\left(0,\s(v)^2\right),
$$
where
\begin{align}
\label{bias}
\b(v)&=\frac{\{1-F_0(v)\}h_{01}''(v)-F_0(v)h_{02}''(v)}{g_1(v)\{1-F_0(v)\}+F_0(v)g_2(v)}\int u^2K(u)\,du\nonumber\\
&\qquad+d_{F_0}(v)\left\{\int_{t=0}^v\frac{\frac{\partial^2}{\partial v^2}h_0(t,v)}{F_0(v)-F_0(t)}\,dt-\int_{u=v}^M\frac{\frac{\partial^2}{\partial v^2}h_0(v,u)}{F_0(u)-F_0(v)}\,du\right\}\int u^2K(u)\,du,
\end{align}
where $h_0$, $h_{01}$ and $h_{02}$ are defined by (\ref{def_h_0}) and
\begin{align}
\label{sigma1}
\s_1(v)=1+d_{F_0}(v)\left\{\int_{t<v}\frac{g(t,v)}{F_0(v)-F_0(t)}\,dt
+\int_{w>v} \frac{g(v,w)}{F_0(w)-F_0(v)}\,dw\right\},
\end{align}
and where $N\left(0,\s(v)^2\right)$ is a normal distribution with first moment zero and variance $\s(v)^2$, defined by
\begin{align}
\label{as-variance}
\s(v)^2=\frac{d_{F_0}(v)}
{\s_1(v)}\int K(u)^2\,du.
\end{align}
\end{theorem}

\begin{remark}
{\rm
Condition (1.1) of \cite{piet:12b} is a separation condition which ensures that the observation intervals $(T_i,U_i)$ do not become arbitrarily small. Because of the additional smoothness conditions, it takes a somewhat more complicated form than the analogous separation condition for the ordinary MLE. We refer for the precise formulation to \cite{piet:12b}.
}
\end{remark}

\section{Deconvolution, smooth functionals}
\label{section:decon}
\setcounter{equation}{0}
The theory of MLEs for deconvolution is full of peculiar facts and unsolved problems. We can again expect that the use of MLEs will produce efficient estimates of smooth functionals. Whether this will give better estimates than naive estimators will depend on the model and in particular on the properties of the tangent spaces, associated with the model. We start with a simple example, where the estimate of the first moment, using the MLE, coincides with a moment estimate.

Suppose our observations $Z_1,\dots,Z_n$ are a sample of the form
$$
Z_i=X_i+Y_i,
$$
where the $X_i$ and $Y_i$ are independent, and $Y_i$ has a (known) normal $N(\m,1)$ distribution. A natural estimate of the first moment of the distribution of the $X_i$ is the estimate
\begin{equation}
\label{moment_estimate}
T_n=n^{-1}\sum_{i=1}^n Z_i-\m.
\end{equation}
The MLE of the unknown distribution function of the $X_i$ is the distribution function $\hat F_n$, maximizing
$$
\ell(F)=\int \log \int\f(z-x-\m)\,dF(x)\,d\H_n(z),
$$
over $F$, where $\H_n$ is the empirical distribution function of the $Z_i$ and $\f$ is the standard normal density. So another estimate of the first moment of the distribution of the $X_i$ is the estimate
$$
T_n'=\int x\,d\hat F_n(x).
$$
But a simple calculation, which is omitted here, shows that, in fact, $T_n'=T_n$, so the two methods produce exactly the same (efficient) estimate here.

This relation does not hold for higher moments however. We could, for example, estimate the variance of the $X_i$ by
$$
U_n=n^{-1}\sum_{i=1}^n \left(Z_i-\bar Z_n\right)^2-1,
$$
where $\bar Z_n$ is the mean of the $Z_i$, but also by
$$
U_n'=\int x^2\,d\hat F_n(x)-\left\{\int x\,d\hat F_n(x)\right\}^2.
$$
Here we do not get $U_n=U_n'$; for example $U_n$ can have negative values, in contrast with $U_n'$.
On theoretical grounds, one would expect the MLE to produce an asymptotically efficient
estimate of the variance, but looking at simulations, one also would expect this
efficiency only to show up for huge sample sizes, because of the highly discrete
character of the MLE, which only has very few points of mass for moderate sample sizes.

The usual method of producing estimates of $F$ is to first estimate the characteristic function of the data in some way, and then use the fact that the characteristic function of the convolution is a product, meaning that one can divide by the characteristic function of the distribution of the known component of the deconvolution to obtain the characteristic function of the unknown component. This does not necessarily produce efficient estimates of the smooth functionals, however, while for the MLE there is a general theory, predicting the efficiency of the estimates of smooth functionals based on the MLE, as also shown in the preceding sections.

Dividing by the characteristic function of the known component becomes more difficult if this characteristic function has zeroes, as in the case of the uniform distribution. In this case the moment estimator (\ref{moment_estimate}) also does not produce an efficient estimate of the first moment. Deconvolution for the case that $Y_i$ has a uniform distribution is sometimes called ``box-car" deconvolution. We consider here the simplest case, where the distributions of the $X_i$ and $Y_i$ both have support $[0,1]$ and the $X_i$ have an absolutely continuous distribution function $F_0$. As discussed in \cite{GrWe:92} (Exercise 2, section 2.3, p.\ 61),
the model is equivalent to the current status model in this case. This is seen in the following way.

Let $\dd_i=1_{\{Z_i\le 1\}}$ and let $Z_i'$ be defined by
\begin{equation}
\label{decon_curstat}
Z_i'=\left\{\begin{array}{lll}
Z_i,\,&\mbox{ if }\dd_i=1,\\
Z_i-1,\,&\mbox{ if }\dd_i=0.
\end{array}
\right.
\end{equation}
Then $Z_1',\dots,Z_n'$ is distributed as a sample from a Uniform$(0,1)$ distribution. Moreover, the log likelihood for the unknown distribution function $F$ can be written
$$
\ell(F)=\sum_{i=1}^n\left\{\dd_i\log F(Z_i')+(1-\dd_i)\log\{1-F(Z_i')\right\},
$$
and we have, for $t,t+h\in(0,1)$
$$
\P\left\{\dd_i=1,\,Z_i'\in [t,t+h]\right\}=\P\left\{Z_i\in [t,t+h]\right\}
\sim h\int_0^t dF_0(u)=hF_0(t),\,h\downarrow0,
$$
and
$$
\P\left\{\dd_i=0,\,Z_i'\in [t,t+h]\right\}=\P\left\{Z_i\in [1+t,1+t+h]\right\}
\sim h\int_t^1 dF_0(u)=h\{1-F_0(t)\},\,h\downarrow0,
$$
so we get factorization of the current status model for $\dd_i$ and the corresponding observation $Z_i'$. This means, by Theorem 5.5 in \cite{GrWe:92}, that
$$
\sqrt{n}\left\{\int x\,d\hat F_n(x)-\int x\,dF_0(x)\right\}\stackrel{{\cal D}}\longrightarrow N(0,\s_{F_0}^2),
$$
where
$$
\s_{F_0}^2=\int_0^1 F_0(t)\{1-F_0(t)\}\,dt.
$$
This also follows from Example 11.2.3b, p. 226, in \cite{sara:00}, where it is at the same time shown that this is the efficient asymptotic variance.

On the other hand, if we would take the moment estimate $T_n$ of (\ref{moment_estimate}) to estimate the first moment of the distribution of the $X_i$, we would get
$$
\sqrt{n}\left\{T_n-\int x\,dF_0(x)\right\}\stackrel{{\cal D}}\longrightarrow N(0,\s^2),
$$
where
$$
\s^2=\frac1{12}+\mbox{var}(X_1).
$$
Since
\begin{align*}
\mbox{var}(X_1)=2\int_0^1 x\{1-F_0(x)\}\,dx-\left\{\int_0^1\{1-F_0(x)\}\,dx\right\}^2,
\end{align*}
a simple variational argument shows that $\s_{F_0}^2<\s^2$, unless $F_0$ is the uniform distribution function, in which case $\s_{F_0}^2=\s^2$. So in this case, the estimate of the first moment, based on the MLE, is more efficient than the moment estimate $T_n$, in contrast with what happened for normal deconvolution.

We now generalize this example to the situation that the convolution kernel is a continuously differentiable decreasing density $g$ on $[0,1]$ and $F_0$ is again an absolutely continuous distribution function, concentrated on $[0,1]$. As in section \ref{section:smooth_functionalsIC}, we consider functionals $K(F)$ satisfying condition (D1), of which the first moment of $F$ is the prototype, so
$$
K(F)=K(F_0)+\int\k_{F_0}(x)\,d\left(F-F_0\right)(x)+O(\|F-F_0\|_2^2),
$$
and we try again to prove:
\begin{equation}
\label{L_2-bound_decon}
\|\hat F_n-F_0\|^2=o_p\left(n^{-1/2}\right),
\end{equation}
and next to prove, using the characterizing properties of the MLE,
\begin{equation}
\label{as_efficiency_relation_decon}
n^{1/2}\int\k_{F_0}(x)\,d\bigl(\hat F_n-F_0\bigr)(x)=n^{1/2}\int\th_{F_0}(z)\,d\left(\H_n-H_0\right)(z)+o_p\left(n^{-1/2}\right).
\end{equation}
Here $H_0$ is the distribution function of the convolution and $\H_n$ the empirical distribution function of the observations $Z_i=X_i+Y_i$, and $\th_{F_0}$ is given by:
\begin{equation}
\label{score_eq}
\th_{F_0}(z)=E\left\{\k_{F_0}(X_1)\bigm|X_1+Y_1=z\right\}.
\end{equation}
Note the analogy with (\ref{as_efficiency_relation1}) and (\ref{canon_gradientIC}) in section \ref{section:smooth_functionalsIC}. Introducing an intermediate function $\f_F$ again, just as in section \ref{section:smooth_functionalsIC}, we get the representation
\begin{equation}
\label{score_decon}
\th_{F_0}(z)=\frac{\int_{u=0}^z g(z-u)\,d\f_{F_0}(u)}{h_0(z)}\,,
\end{equation}
where $\f_{F_0}$ is absolutely continuous w.r.t. $F_0$, $g$ is the (decreasing) density of the $Y_i$ on $[0,1]$, and $h_0$ is the density of the observations $Z_i$. This leads to the integral equation (in $\f_{F_0}$):
\begin{equation}
\label{adjoint_eq}
\k_{F_0}(x)=E\left\{\th_{F_0}(Z_1)|X_1=x\right\}=\int_{z\ge x}\frac{\int_{u=0}^z g(z-u)\,d\f_{F_0}(u)}{h_0(z)}\,g(z-x)\,dz,
\end{equation}
which is the adjoint equation to the equation (\ref{score_eq}) (more on the relation between the equation (\ref{score_eq}) and its adjoint (\ref{adjoint_eq}) in the deconvolution problem can be found in part 1 of \cite{GrWe:92}). Note that this is an equation in the function $\f_{F_0}'$ (or the measure $d\f_{F_0}$) rather than the function $\f_{F_0}$. Also note that we can write
\begin{equation}
\label{adjoint_eq2}
\k_{F_0}(x)=\int_{z\ge x}\th_{F_0}(z)\,g(z-x)\,dz,
\end{equation}
giving a seemingly simpler equation in $\th_{F_0}$. The essential (sometimes ignored) fact, however, is that $\th_{F_0}$ has to have a representation of the form (\ref{score_decon}) (or, more generally, that it has to be a limit of representations of this form: it has to be in the closure of the range of the score operator). Without this restriction on $\th_{F_0}$ the equation (\ref{adjoint_eq2}) would have infinitely many solutions; the restriction that $\th_{F_0}$ has to have a representation of the form (\ref{score_decon}) or has to be a limit of such representations makes the solution unique, however.

Defining $\f_{F_0}(0)=0$ and using integration by parts we get:
\begin{align*}
\frac{\int_{u=0}^z g(z-u)\,d\f_{F_0}(u)}{h_0(z)}
=\frac{g(0)\f_{F_0}(z)+\int_0^z \f_{F_0}(u)g'(z-u)\,du}{h_0(z)}\,
\end{align*}
and by differentiating (\ref{adjoint_eq}) we obtain the integral equation
\begin{equation}
\label{decon_inteq}
\f_{F_0}(x)+\int_{u=0}^1 A(x,u)\f_{F_0}(u)\,du=-\frac{h_0(x)\k_{F_0}'(x)}{g(0)^2},\,x\in(0,1),
\end{equation}
where the kernel $A$ of the integral equation is given by:
\begin{equation}
\label{def_A}
A(x,u)=\frac{g'(x-u)}{g(0)}+\frac{h_0(x)g'(u-x)}{g(0)h_0(u)}+\frac{h_0(x)}{g(0)^2}\int_{z=x\vee u}^{1+x\wedge u}\frac{g'(z-u)g'(z-x)}{h_0(z)}\,dz.
\end{equation}
So, just as in section \ref{section:smooth_functionalsIC}, we obtain a Fredholm integral equation of the second kind.

To give a concrete example of what these functions $\th_{F_0}$ and $\f_{F_0}$ look like, we take $g$ equal to the ``elbow density"
$g(x)=2(1-x)1_{[0,1]}(x)$ and let $F_0$ be the uniform distribution function on $[0,1]$. Furthermore, we take
\begin{equation}
\label{elbow_gradient}
\k_{F_0}(x)=x-\int u\,dF_0(u),
\end{equation}
the gradient corresponding to the first moment of the distribution, given by $F_0$. For this model the kernel $A(x,u)$ becomes:
$$
A(x,u)=-1_{[u,1)}(x)-\frac{h_0(x)1_{[x,1)}(u)}{h_0(u)}+h_0(x)\int_{z=x\vee u}^{1+x\wedge u}\frac{dz}{h_0(z)}\,.
$$
Moreover:
\begin{equation}
\label{h_0-elbow1}
h_0(z)=\left\{\begin{array}{lll}
z(2-z)&,\,z\in[0,1],\\
(2-z)^2&,\,z\in(1,2],
\end{array}
\right.
\end{equation}
So we get:
\begin{align}
\label{kernel_A}
A(x,u)=-1_{(0,x)}(u)-\frac{x(2-x)1_{[x,1)}(u)}{u(2-u)}+x(2-x)\left\{\tfrac12\log\left(\frac{2-x\vee u}{x\vee u}\right)+\frac{x\wedge u}{1-x\wedge u}\right\},
\end{align}
for $x,u\in(0,1)$. This immediately points to one reason why solving this type of integral equation is more difficult than solving the integral equation we studied in section \ref{section:smooth_functionalsIC}: the kernel of the integral equation is unbounded.

\begin{figure}[!ht]
\begin{center}
\includegraphics[scale=0.5]{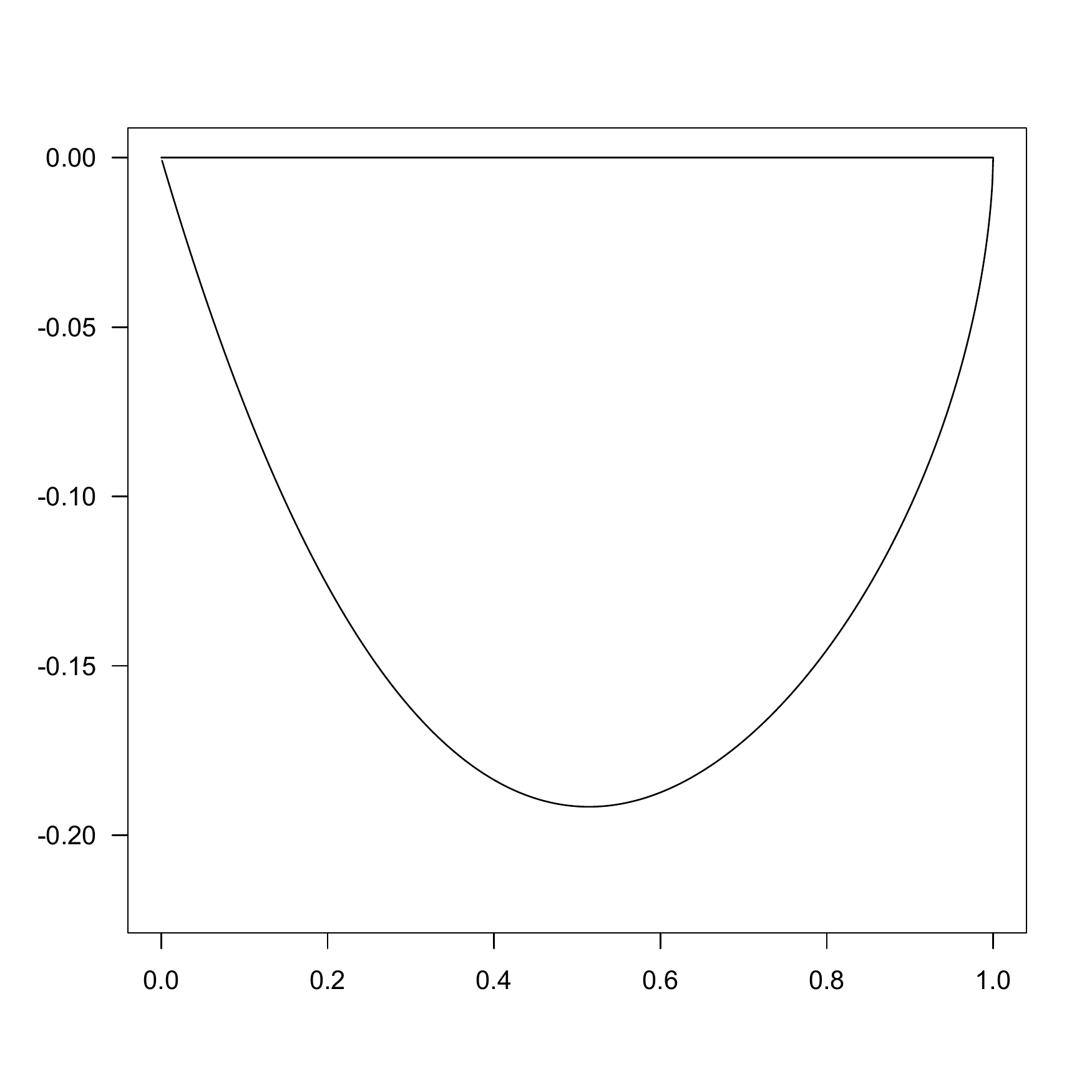}
\end{center}
\caption{The function $\f_{F_0}$, solving the integral equation (\ref{decon_inteq2}).}
\label{fig:phi_decon}
\end{figure}

\begin{figure}[!ht]
\begin{center}
\includegraphics[scale=0.5]{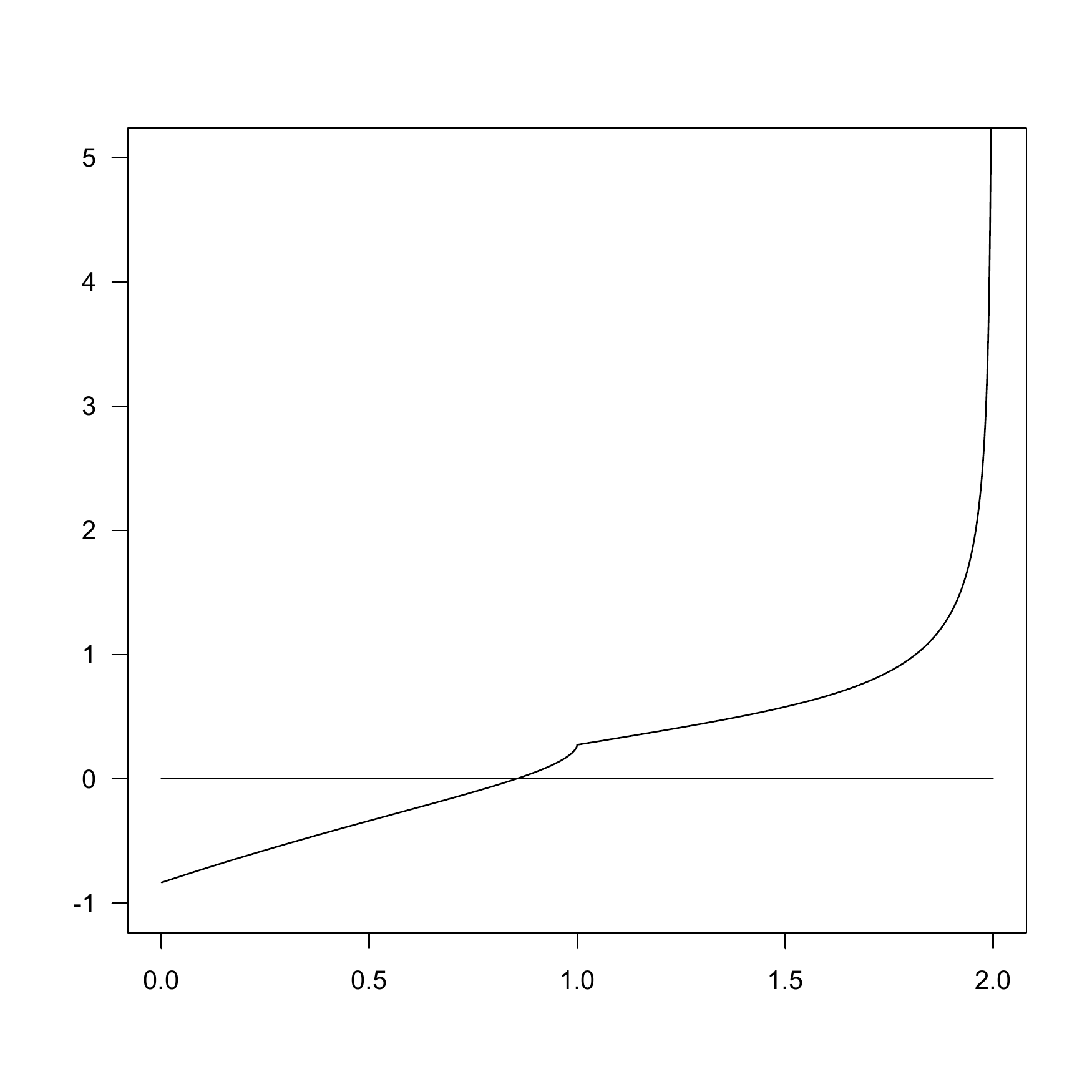}
\end{center}
\caption{The function $\th_{F_0}$, defined by (\ref{score_decon}), for the same model as used in Figure \ref{fig:phi_decon}.}
\label{fig:theta_decon}
\end{figure}

\begin{figure}[!ht]
\begin{center}
\includegraphics[scale=0.5]{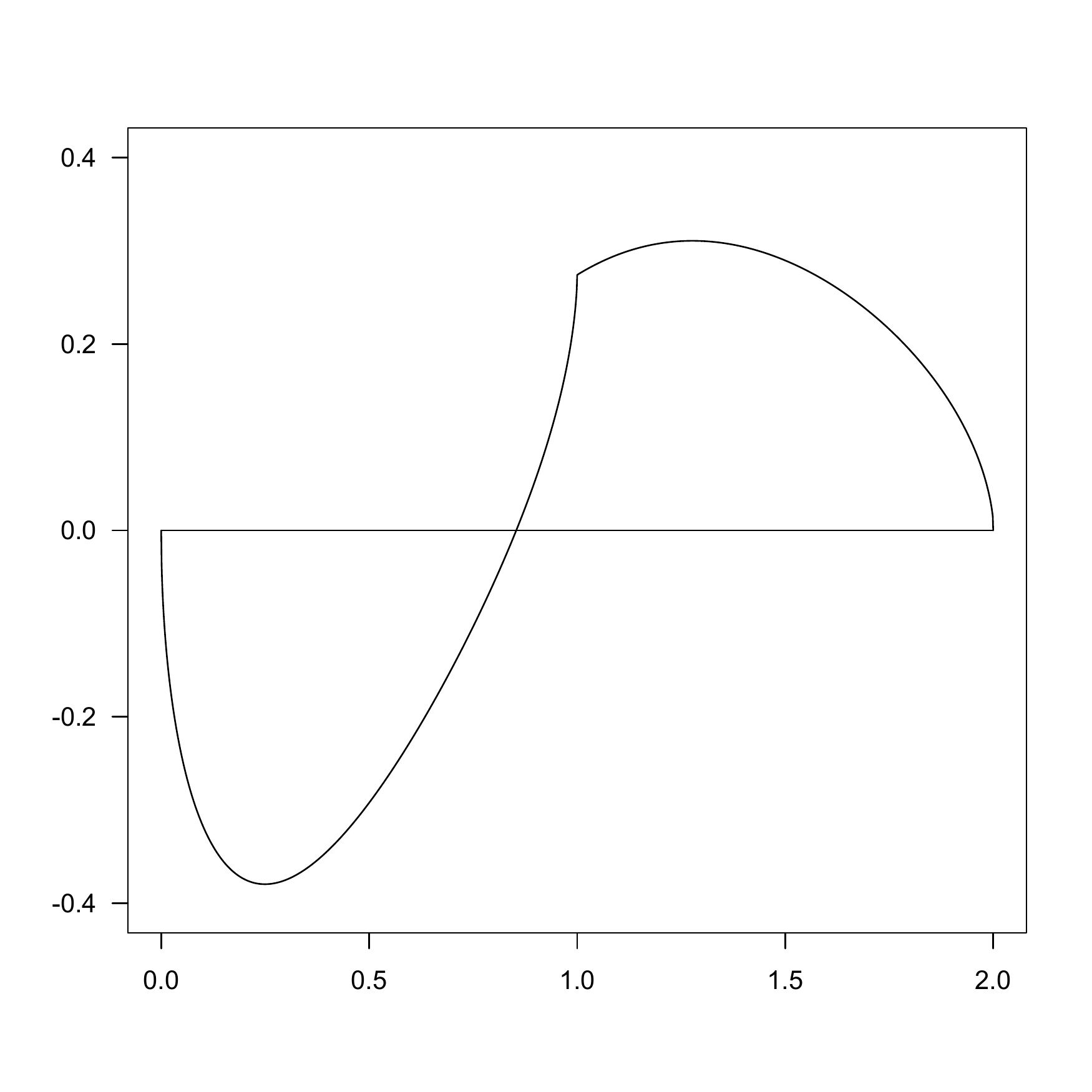}
\end{center}
\caption{The function $\th_{F_0}\sqrt{h_0}$, for the same model as used in Figure \ref{fig:phi_decon}.}
\label{fig:theta_sqrt(h)}
\end{figure}

Furthermore, taking the gradient given by (\ref{elbow_gradient}) and using (\ref{h_0-elbow1}) again, the integral equation (\ref{decon_inteq}) turns into
\begin{equation}
\label{decon_inteq2}
\f_{F_0}(x)+\int_{u=0}^1 A(x,u)\f_{F_0}(u)\,du=-\tfrac14x(2-x),\,x\in(0,1).
\end{equation}
This equation was solved numerically and a picture of $\f_{F_0}$ is given in Figure \ref{fig:phi_decon}. The corresponding $\th_{F_0}$, defined by (\ref{score_decon}), is shown in Figure \ref{fig:theta_decon}. The function $z\mapsto \th_{F_0}(z)$ is unbounded near $z=2$ and has a cusp at $z=1$. 
It can be shown, though, that
\begin{equation}
\label{eff_var_elbow}
\int\th_{F_0}(z)^2 h_0(z)\,dz<\infty.
\end{equation}
and $\th_{F_0}\in L_2^0(H_0)$, where $L_2^0(H_0)$ is the space of square integrable function w.r.t.\ $dH_0$ which integrate to zero w.r.t. $dH_0$. In fact the left-hand side of (\ref{eff_var_elbow}) can be expected to be the (efficient) asymptotic variance of
$$
n^{1/2}\left\{\int x\,d\hat F_n(x)-\int x\,dF_0(x)\right\}.
$$
A picture of $\th_{F_0}\sqrt{h_0}$ is shown in Figure \ref{fig:theta_sqrt(h)}. Numerical evaluation of (\ref{eff_var_elbow}) gave indeed a value of approximately $0.137$, and, for example, a simulation of $1000$ samples of size $n=1000$, where the MLE was computed using the support reduction algorithm of \cite{piet_geurt_jon:08}, gave the value $0.139$ for $n$ times the variance of the sample estimates (we have the impression that the sample variance times $n$ converge to the asymptotic value from above, as $n\to\infty$; sample size $n=100$ gave slightly larger values), so simulations seem to give a nice agreement with the conjectured asymptotic variance.

In order to prove
$$
n^{1/2}\left\{\int x\,d\hat F_n(x)-\int x\,dF_0(x)\right\}\stackrel{\cal D}\longrightarrow
N(0,\s_0^2),
$$
where
$$
\s_0^2=\int\th_{F_0}(z)^2 h_0(z)\,dz,
$$

we will need the following lemma, which is analogous to Lemma \ref{lemthF}.

\begin{lemma}
\label{lemthF2}
Let the density $\hat h_n$ be defined by
$$
\hat h_n(z)=\int g(z-x)\,d{\hat F_n}(x),\,z\in[0,2].
$$
and let $\t_1$ and $\t_m$ be the smallest and largest point of jump of $\hat F_n$, respectively.
Furthermore, let, in analogy with (\ref{score_decon}), $\th_{\hat F_n}$ be defined by
\begin{equation}
\label{thetaFn_representation_MLE_decon}
\th_{\hat F_n}(z)=\frac{\int_{u=0}^z g(z-u)\,d\f_{\hat F_n}(u)}{\hat h_n(z)}\,,\,z\in[\t_1,1+\t_m),
\end{equation}
where the function $\f_{\hat F_n}$ solves the integral equation (in $\f$)
\begin{equation}
\label{decon_inteq3}
\f(x)+\int_{u=\t_1}^{\t_m} A_n(x,u)\f(u)\,du=\frac{\hat h_n(x)}{g(0)^2},\,x\in[\t_1,\t_m),
\end{equation}
and where the kernel $A_n$ is given by
$$
A_n(x,u)=\frac{g'(x-u)}{g(0)}+\frac{\hat h_n(x)g'(u-x)}{g(0)\hat h_n(u)}+\frac{\hat h_n(x)}{g(0)^2}\int_{z=x\vee u}^{1+x\wedge u}\frac{g'(z-u)g'(z-x)}{\hat h_n(z)}\,dz.
$$
We define $\f_{\hat F_n}(x)=0$, if $x\in[0,\t_1)$ or $x\ge\t_m$. Then $\th_{\hat F_n}$ satisfies
$$
\int_{z\ge x}\theta_{\hat F_n}(z)g(z-x)\,dz=x-\int u\,d\hat F_n(u),\,x\in[\t_1,\t_m).
$$
The function $\th_{\hat F_n}$ can uniquely and continuously be extended to $(0,2\vee(1+\t_m))$ such that
\begin{equation}
\label{theta_extension}
\int_{z\ge x}\theta_{\hat F_n}(z)g(z-x)\,dz=x-\int u\,d\hat F_n(u),\,x\in(0,1).
\end{equation}
Moreover, for $\th_{\hat F_n}$, extended in this way, we have:
\begin{equation}
\label{obs_space_representation_decon}
\int x\,d\bigl(\hat F_n-F_0\bigr)(x)=-\int \theta_{\hat F_n}(z)\,dH_0(z).
\end{equation}
\end{lemma}

\begin{remark}
{\rm Note that
$$
\hat h_n(z)=0,\,z\notin[\t_1,1+\t_m).
$$
Also note that we can have $\t_m<1$ and $\t_m>1$, the latter event is more typical.
}
\end{remark}

\vspace{0.3cm}
To show that (\ref{obs_space_representation_decon}) holds, note that
\begin{align*}
&\int \theta_{\hat F_n}(z)\,dH_0(z)=\int \theta_{\hat F_n}(z)\int g(z-x)\,dF_0(x)\,dz\\
&=\int\left\{\int\theta_{\hat F_n}(z)g(z-x)\,dz\right\}\,dF_0(x)
=\int\left\{x-\int x\,d\hat F_n(x)\right\}\,dF_0(x)=-\int x\,d\bigl(\hat F_n-F_0\bigr)(x).
\end{align*}

\begin{figure}[!ht]
\begin{center}
\includegraphics[scale=0.5]{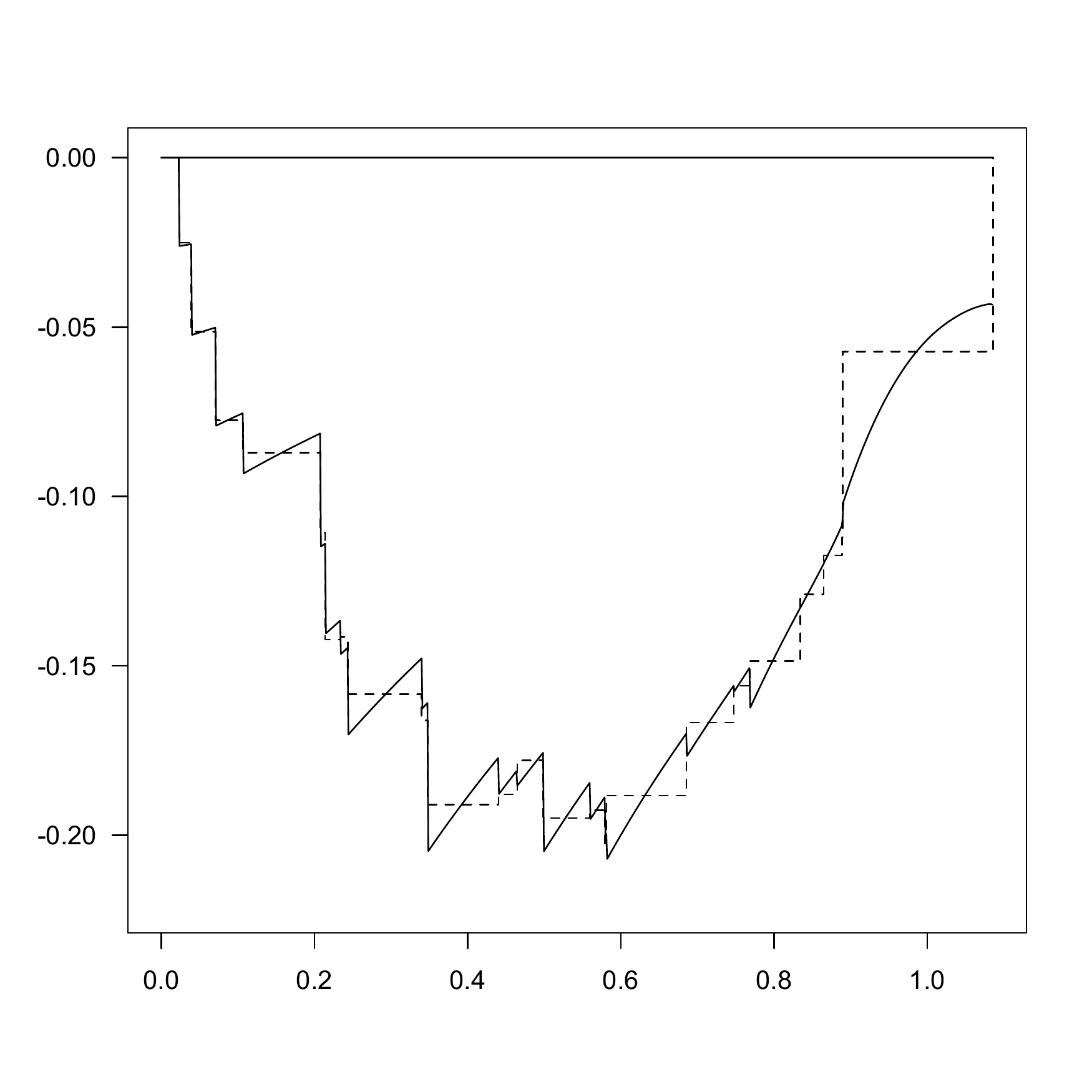}
\end{center}
\caption{The function $\f_{\hat F_n}$ (solid), defined as in Lemma \ref{lemthF2}, and the function $\bar\f_{\hat F_n}$ (dashed), defined as the solution of the equation (\ref{adjoint_eq3}), for the same model as used in Figure \ref{fig:phi_decon}, and for a sample of size $n=1000$; $\t_m\approx1.08617$.}
\label{fig:phi_n_decon}
\end{figure}

Just as in section \ref{section:smooth_functionalsIC}, we would like to have a relation of the form
$$
\int x\,d\bigl(\hat F_n-F_0\bigr)(x)=\int\theta_{\hat F_n}(z)\,d\bigl(\H_n-H_0\bigr)(z)+o_p\left(n^{-1/2}\right)
$$
instead of (\ref{obs_space_representation_decon}). Proceeding as in section \ref{section:smooth_functionalsIC}, we construct a function $\bar\f_{\hat F_n}$ which is absolutely continuous w.r.t.\ $\hat F_n$ and which is close to $\f_{\hat F_n}$. Next we define
\begin{equation}
\label{def_bar_theta}
\bar\th_{\hat F_n}(z)=\left\{\begin{array}{lll}
\displaystyle{\frac{\int_{x\in[0,z]}g(z-x)\,d\bar\f_{\hat F_n}(x)}{\hat h_n(z)}}\,&,\,z\in[\t_1,1+\t_m),\\
0\,&,\mbox{ otherwise.}
\end{array}
\right.
\end{equation}
To this end, we define $\bar\f_{\hat F_n}$ as the solution of the equation in $\f$
\begin{equation}
\label{adjoint_eq3}
\int_{u\in[\t_1,\t_m]}\left\{\int_{z\ge x\vee u}\frac{g(z-u)g(z-x)}{\hat h_n(z)}\,dz\right\}\,d\f(u)=x-\int x\,d\hat F_n(x)\qquad\,,\mbox{a.e.\ }\left[d\hat F_n\right].
\end{equation}
Note that this is a finite matrix equation which only has to be solved at the points of jump of $\hat F_n$.
Then $\bar\f_{\hat F_n}$ is absolutely continuous w.r.t.\ $\hat F_n$, since it is a right-continuous piecewise constant function, having (finitely many) jumps at the same locations as $\hat F_n$. 
It follows from Proposition 2.1, p.\ 54 in \cite{GrWe:92} (see also \cite{steffi:09b} and \cite{steffi:12}) that
$$
\int_{z\ge\t_i}\frac{g(z-\t_i)}{\hat h_n(z)}\,d\H_n(z)=1,\,i=1,\dots,m.
$$
So we get:
\begin{align}
\label{MLE_relation_decon}
\int\bar\th_{\hat F_n}(z)\,d\H_n(z)&=\int\left\{\int \frac{g(z-x)}{\hat h_n(z)}\,d\H_n(z)\right\}\,d\bar\f_{\hat F_n}(x)\nonumber\\
&=\int_{x\in[\t_1,\t_m]}\,d\bar\f_{\hat F_n}(x)=\bar\f_{\hat F_n}(\t_m)-\bar\f_{\hat F_n}(0)=0.
\end{align}
A picture of the functions $\f_{\hat F_n}$ and $\bar\f_{\hat F_n}$ is shown in Figure \ref{fig:phi_n_decon}.

Furthermore, by Theorem 5, p.\ 522, of \cite{steffi:09}:
$$
\|\hat F_n-F_0\|_2=O_p\left(n^{-1/3}\right).
$$
This suggests, using methods analogous to the methods used in \cite{GeGr:97},
\begin{equation}
\label{L_2-bound_dec}
\int\bigl\{\bar\th_{\hat F_n}(z)-\th_{\hat F_n}(z)\bigr\}\,dH_0(z)=o_p\left(n^{-1/2}\right).
\end{equation}
In fact, we have:
\begin{align*}
&\int\bigl\{\bar\th_{\hat F_n}(z)-\th_{\hat F_n}(z)\bigr\}\hat h_n(z)\,dz=\int_{\t_1}^{1+\t_m}\bigl\{\bar\th_{\hat F_n}(z)-\th_{\hat F_n}(z)\bigr\}\hat h_n(z)\,dz\\
&=\int_{x\in[\t_1,\t_m]}\left\{\int_{z\ge x} g(z-x)\,dz\right\}\,d\left(\bar\f_{\hat F_n}-\f_{\hat F_n}\right)(x)=\int_{x\in[\t_1,\t_m]}\,d\left(\bar\f_{\hat F_n}-\f_{\hat F_n}\right)(x)=0,
\end{align*}
and hence
\begin{align*}
&\int\bigl\{\bar\th_{\hat F_n}(z)-\th_{\hat F_n}(z)\bigr\}h_0(z)\,dz
=\int\bigl\{\bar\th_{\hat F_n}(z)-\th_{\hat F_n}(z)\bigr\}\bigl\{h_0(z)-\hat h_n(z)\bigr\}\,dz.
\end{align*}
We write the integrand on the right-hand side as the product of the functions
$$
z\mapsto \left\{\bar\th_{\hat F_n}(z)-\th_{\hat F_n}(z)\right\}\left\{\sqrt{\hat h_n(z)}+\sqrt{h_0(z)}\right\}
$$
and
$$
z\mapsto \sqrt{\hat h_n(z)}-\sqrt{h_0(z)}
$$

\begin{figure}[!ht]
\begin{center}
\includegraphics[scale=0.5]{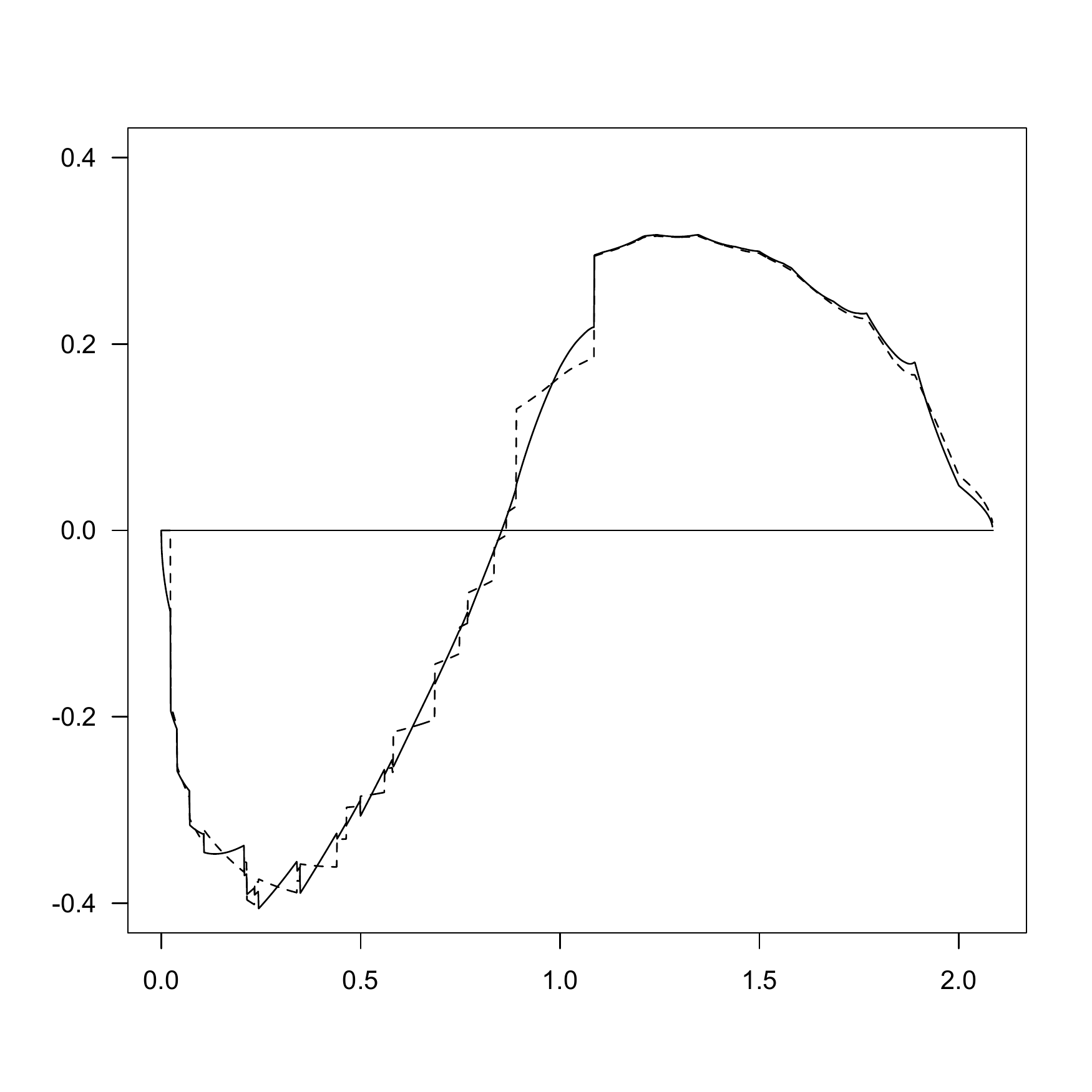}
\end{center}
\caption{The function $\th_{\hat F_n}\{\sqrt{\hat h_n}+\sqrt{h_0}\}/2$ (solid), defined by (\ref{thetaFn_representation_MLE_decon}),  and the function $\bar\th_{\hat F_n}\{\sqrt{\hat h_n}+\sqrt{h_0}\}/2$ (dashed), defined by (\ref{def_bar_theta}), for the same model as used in Figure \ref{fig:phi_decon}, and for the same sample as used in Figure \ref{fig:phi_n_decon}; $\t_m\approx1.08617$.}
\label{fig:theta_n_decon}
\end{figure}

\begin{figure}[!ht]
\begin{center}
\includegraphics[scale=0.5]{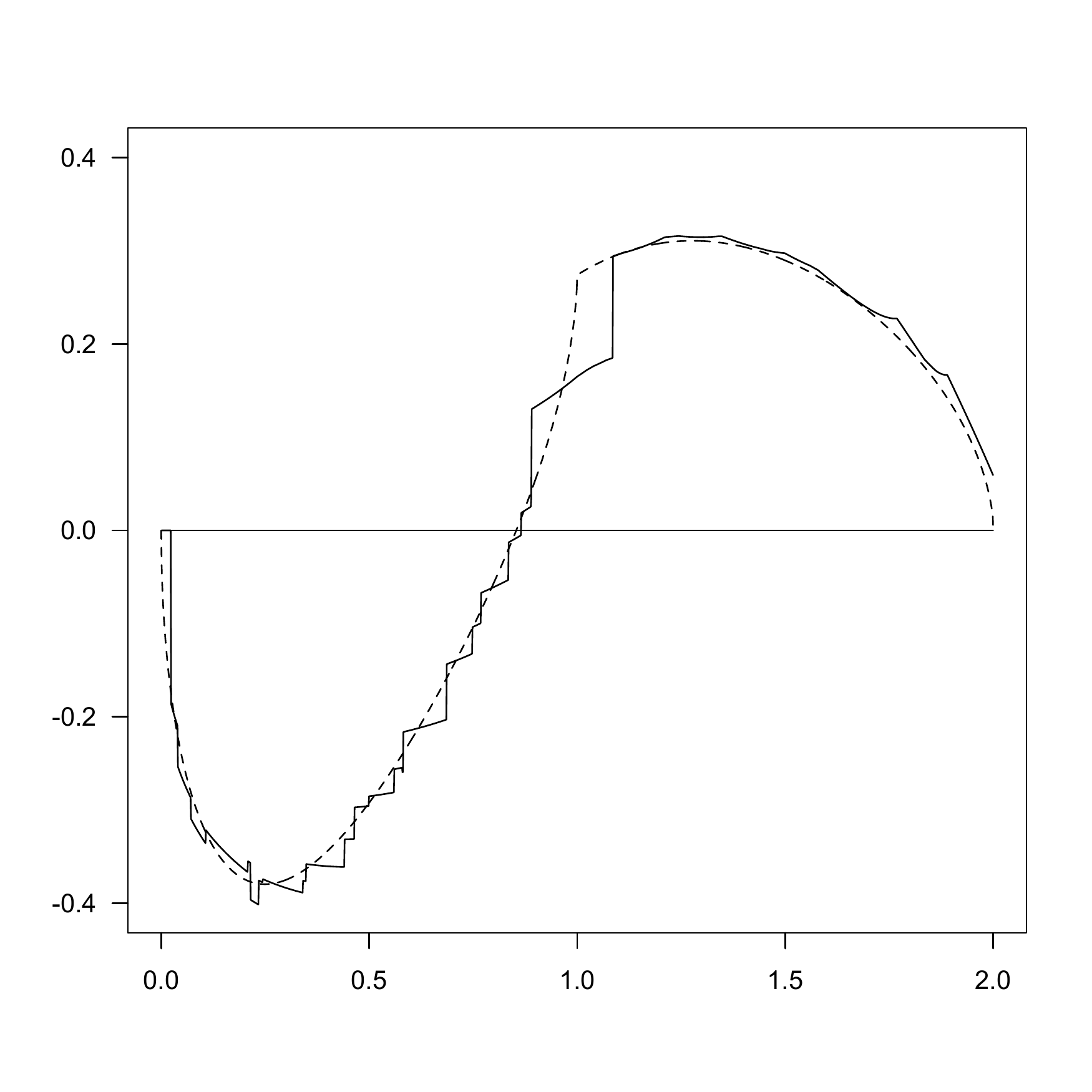}
\end{center}
\caption{The function $\bar\th_{\hat F_n}\{\sqrt{\hat h_n}+\sqrt{h_0}\}/2$ (solid), defined by (\ref{thetaFn_representation_MLE_decon}),  and the function $\th_{F_0}\sqrt{h_0}$ (dashed), defined by (\ref{score_decon}), on the interval $[0,2]$, for the same model as used in Figure \ref{fig:phi_decon}, and for the same sample as used in Figure \ref{fig:phi_n_decon}.}
\label{fig:bar_theta_theta0}
\end{figure}

It is proved in \cite{steffi:09} that
$$
\left\{\int\left\{\sqrt{\hat h_n(z)}-\sqrt{h_0(z)}\right\}^2\,dz\right\}^{1/2}=O_p\left(n^{-1/3}\right).
$$
If it can also be shown that
$$
\left\{\int \left\{\bar\th_{\hat F_n}(z)-\th_{\hat F_n}(z)\right\}^2\left\{\hat h_n(z)+h_0(z)\right\}\,dz\right\}^{1/2}=o_p\left(n^{-1/6}\right),
$$
we would obtain, using (\ref{obs_space_representation_decon}) to (\ref{L_2-bound_dec}) and the Cauchy-Schwarz inequality,
\begin{align}
\label{decon_rep_obs}
\int x\,d\bigl(\hat F_n-F_0\bigr)(x)&=-\int \theta_{\hat F_n}(z)\,dH_0(z)
=-\int \bar\theta_{\hat F_n}(z)\,dH_0(z)+o_p\left(n^{-1/2}\right)\nonumber\\
&=\int \bar\theta_{\hat F_n}(z)\,d\bigl(\H_n-H_0\bigr)(z)+o_p\left(n^{-1/2}\right).
\end{align}
A picture of the functions $\th_{\hat F_n}\{\sqrt{\hat h_n}+\sqrt{h_0}\}/2$ and $\bar\th_{\hat F_n}\{\sqrt{\hat h_n}+\sqrt{h_0}\}/2$ is shown in Figure \ref{fig:theta_n_decon}, showing that these function are really close on the interval $[0,1+1\vee\t_m]$. Figure \ref{fig:bar_theta_theta0} compares $\bar\th_{\hat F_n}\{\sqrt{\hat h_n}+\sqrt{h_0}\}/2$ and $\th_{F_0}\sqrt{h_0}$.

The only remaining step would be to show that
$$
\sqrt{n}\int \bar\theta_{\hat F_n}(z)\,d\bigl(\H_n-H_0\bigr)(z)
=\sqrt{n}\int \theta_{F_0}(z)\,d\bigl(\H_n-H_0\bigr)(z)+o_p(1),
$$
and this representation would again give asymptotic efficiency of the estimate of the first moment, based on the MLE.
It is clear that the heart of the difficulty of the proof is the unboundedness of the functions $\th_{F_0}$, $\th_{\hat F_n}$ and $\bar\th_{\hat F_n}$ at the right endpoint of the interval on which they are defined, which is caused by the fact that the decreasing convolution density $g$ approaches zero at the right endpoint of the interval $[0,1]$.

\section{Deconvolution, local limits}
\label{section:decon_local}
\setcounter{equation}{0}
We briefly discuss the conjectured local limit behavior of the MLE. Consistency of the MLE for deconvolution has been proved quite generally in \cite{steffi:12}. For the local limit distribution theory for the case that the mixture density is decreasing the following conjecture was launched.

\begin{conjecture}
\label{th:conjecture}
\mbox{\rm (Conjectured theorem  in 5.4 \cite{Gr:91}), also Theorem  in 5.4 in \cite{GrWe:92}.)} Let $g$ be a right-continuous decreasing density on
$[0,\infty)$, having only a finite number of discontinuity points
$a_0=0<a_1<\dots<a_m$. Moreover, suppose that
$g$ has a derivative $g'(x)$ at points $x\ne a_i,\,i=0,\dots,m$, satisfying
$$
\int_{(0,\infty)}\frac{g'(x)^2}{g(x)}\,dx<\infty,
$$
where the integrand is defined to be zero at the points $a_i$ and at points $x$
where $g$ is zero, and where $g'$ is bounded and
continuous on the intervals $(a_{i-1},a_i)$, $i=1,\dots,m+1$, with
$a_{m+1}\stackrel{{\rm def}}{=}\infty$.

Furthermore, assume that there exist positive constants  $k_1$ and $k_2$ such
that the derivative $g'$ of $g$ satisfies the relation
$$
|g'(t+u)|\le k_1|g'(t)|,
$$
for all $t>0$ and $0<u<k_2$, such that $a_i<t<t+u<a_{i+1}$ for some $i$,
$0\le i\le m$. 

Let the convolution density $h$ be given by
$$
h(z)=\int g(z-x)\,dF_0(x),\,z\ge 0,
$$
where the distribution function $F_0$ of the (non-negative) random variables
$X_i,1\le i\le n$, is continuously differentiable at
$z_0>0$, with derivative $f_0(z_0)>0$ at $z_0$. Then
$$
n^{1/3}\left\{F_n^{(1)}(z_0)-F_0(z_0)\right\}
f_0(z_0)^{-1/3}\Bigl\{2\sum_{i=0}^m
\bigl(g(a_i)-g(a_i-)\bigr)^2\big/h(z_0+a_i)\Bigr\}^{1/3}
\stackrel{\cal D}{\longrightarrow}2Z,
$$
where $\stackrel{\cal D}{\longrightarrow}$ denotes convergence in distribution, and where
$Z$ is the last time where standard two-sided Brownian motion minus the parabola
$y(t)=t^2$ reaches its maximum.
\end{conjecture}

Specializing $g$ to the Uniform distribution on $[0,1]$, the discontinuity points are $a_0=0$ and $a_1=1$, and we get for a fixed interior point $t_0\in(0,1)$:
\begin{equation}
\label{local_limit_unif_decon}
n^{1/3}\bigl\{\hat F_n(t_0)-F_0(t_0)\bigr\}\big/
\bigl\{\tfrac12F_0(t_0)\left(1-F_0(t_0)\right)f_0(t_0)\bigr\}^{1/3}
\stackrel{\cal D}{\longrightarrow}2Z,
\end{equation}
where $Z$ is as in the conjectured Theorem \ref{th:conjecture}. We know this to be true by the interpretation in terms of the current status model, see (\ref{decon_curstat}) in section \ref{section:decon}.

Specializing $g$ to the standard exponential distribution on $[0,\infty)$, gives only one discontinuity point $a_0=0$, and the conjectured theorem then yields, for $t_0>0$:
\begin{equation}
\label{local_limit_unif_expon}
n^{1/3}\bigl\{\hat F_n(t_0)-F_0(t_0)\bigr\}\big/
\bigl\{\tfrac12f_0(t_0)h(t_0)\bigr\}^{1/3}\stackrel{\cal D}{\longrightarrow}2Z.
\end{equation}
This is also proved to be true in \cite{geurt:98b}.

For the more general case of a decreasing density, the conjectured theorem still has not been proved. The problem is discussed in some detail in Chapter 5 of \cite{steffi:09b}. Possibly the assumptions are somewhat too strong. In Chapter 5 of \cite{steffi:09b} the following conditions are used.
\begin{enumerate}
\item[(i)] The distribution function $F_0$ is continuous on $[0,S_0]$, and $F_0(S_0)=1$, where $S_0<\infty$.
\item[(ii)] In a neighborhood of $t_0\in(0,S_0)$, $F_0$ is continuously differentiable, with derivative $f_0$, satisfying
$$
f_0(t_0)>0.
$$
\item[(iii)] The density $g$ is bounded, decreasing and continuous on $[0,\infty)$ and has compact support $[0,S_g]$.
Moreover, $g$ has a bounded Lipschitz continuous derivative on $(0,S_g)$.
\end{enumerate}

Under these assumptions, proving the conjectured Theorem \ref{th:conjecture} again depends on being able to deal with certain integral equations. Just as in the case of interval censoring, Theorem \ref{th:conjecture} will follow if we can prove that a certain remainder term is a smooth functional, which converges to zero at a faster rate than the cube root $n$ rate. The analysis starts by considering the estimate $\hat h_n$ of the density $h_0$, generating the observations $Z_i=X_i+Y_i$, 
$$
\hat h_n(t)=\int g(t-x)\,d\hat F_n(x).
$$
where $\hat F_n$ is the MLE. We write this in the form
\begin{equation}
\label{int_by_parts}
\hat h_n(t)=g(0)\hat F_n(t)-\int_{x=0}^t\{g(t-x)-g(0)\}\,d\hat F_n(x).
\end{equation}
Note that the second term on the right-hand side of (\ref{int_by_parts}) has a {\it continuous} (but not differentiable) integrand at $t$, in contrast with the preceding representation of $\hat h_n$. We expect
\begin{equation}
\label{smooth_local}
\int_{x=0}^t\{g(t-x)-g(0)\}\,d\hat F_n(x)-\int_{x=0}^t\{g(t-x)-g(0)\}\,dF_0(x)=O_p\left(n^{-1/2}\right),
\end{equation}
to be a smooth functional, which would mean that the dominant asymptotic behavior of $\hat h_n(t)$ is given by
\begin{equation}
\label{dominant_behavior}
g(0)\hat F_n(t)+\int_{x=0}^t\{g(t-x)-g(0)\}\,dF_0(x).
\end{equation}
This, in turn, would imply that the MLE would be asymptotically equivalent to a certain ``toy estimator", obtained by doing one step of the iterative convex minorant algorithm (for a discussion of this algorithm, see \cite{geurt:98}), and then the conjectured theorem would hold, as further explained in \cite{Gr:91} and \cite{GrWe:92}. This would give for the present model:
\begin{equation}
\label{local_limit_main}
n^{1/3}\left(\frac{2g(0)^2}{f_0(t_0)h_0(t_0)}\right)^{1/3}\bigl\{\hat F_n(t_0)-F_0(t_0)\bigr\}
\stackrel{\cal D}{\longrightarrow}2Z,
\end{equation}
where $Z$ is as in the conjectured Theorem \ref{th:conjecture}.

\begin{figure}[!ht]
\begin{center}
\includegraphics[scale=0.5]{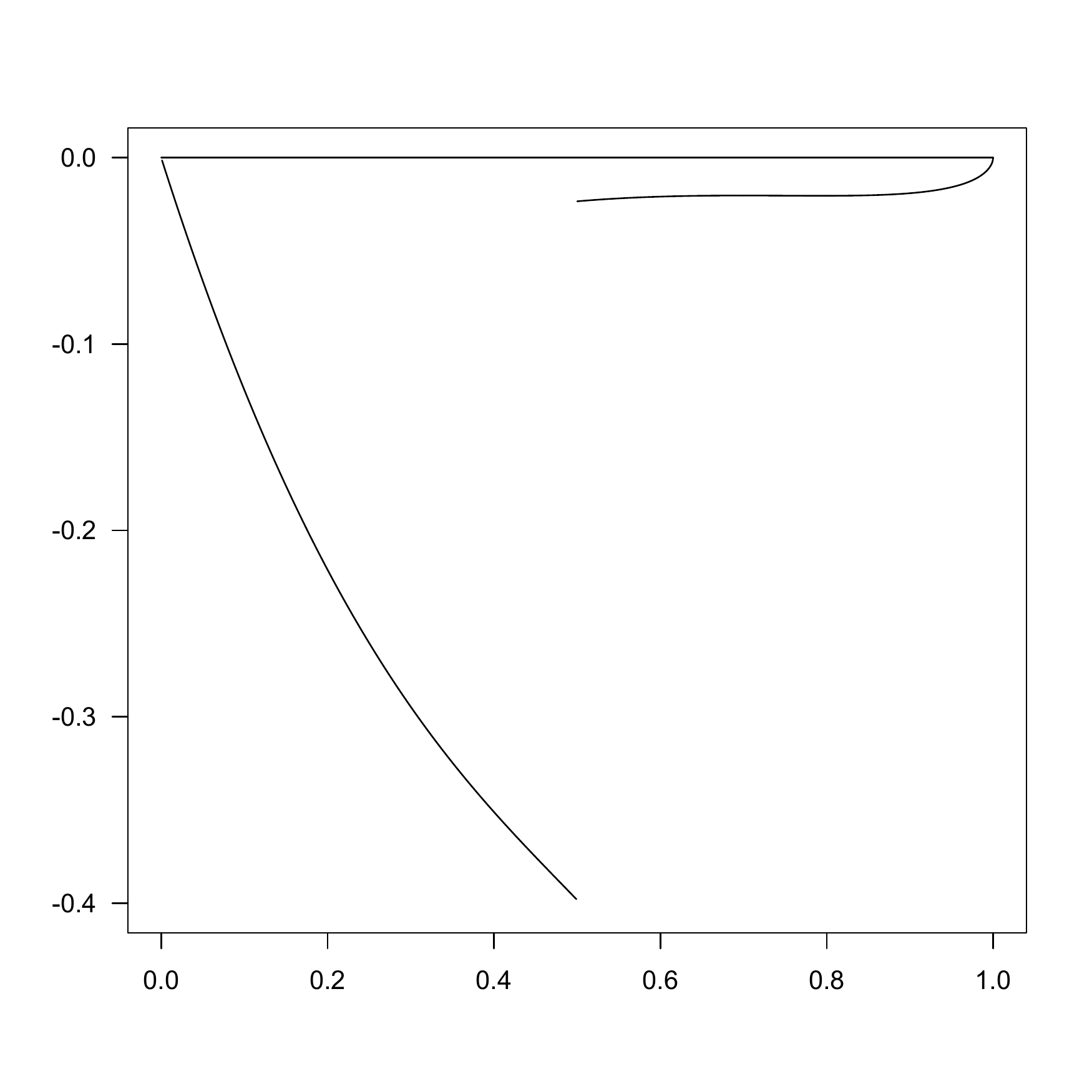}
\end{center}
\caption{The function $\f_{t,F_0}$, solving the integral equation (\ref{Fredholm_local}) for $t=0.5$, where $F_0$ is the uniform distribution function on $[0,1]$ and $g$ the elbow density $g(x)=2(1-x)1_{[0,1]}(x)$.}
\label{fig:phi_local}
\end{figure}

To prove (\ref{smooth_local}), we this time have to deal with the integral equation (\ref{adjoint_eq}) in $\f$:
$$
\k_{t,F_0}(x)=\int_{z\ge x}\frac{\int_{u=0}^z g(z-u)\,d\f(u)}{h_0(z)}\,g(z-x)\,dz,
$$
where in this case:
\begin{equation}
\label{adjoint_eq4}
\k_{F_0}(t,x)=\{g(t-x)-g(0)\}1_{[0,t)}(x)-\int_{u=0}^t \{g(z-u)-g(0)\}\,dF_0(u).
\end{equation}
Assuming $S_0=S_g=1$ in the conditions above, we get that differentiation leads again to a Fredholm integral equation of the second kind:
\begin{equation}
\label{Fredholm_local}
\f_{t,F_0}(x)+\int_{u=0}^1 A(x,u)\f_{t,F_0}(u)\,du=\frac{h_0(x)g'(t-x)1_{[0,t)}(x)}{g(0)^2},\,x\in(0,1),
\end{equation}
where the kernel $A$ is given by (\ref{def_A}). But note that this time the expression on the right-hand side has a jump discontinuity at $t$. This also leads to a solution with a jump discontinuity at the same location.
For the example where $F_0$ is the uniform distribution function on $[0,1]$ and $g$ the ``elbow density" $g(x)=2(1-x)1_{[0,1]}(x)$ the solution is given in Figure \ref{fig:phi_local}, taking $t=0.5$. Note the jump at $t=0.5$.

\begin{figure}[!ht]
\begin{center}
\includegraphics[scale=0.5]{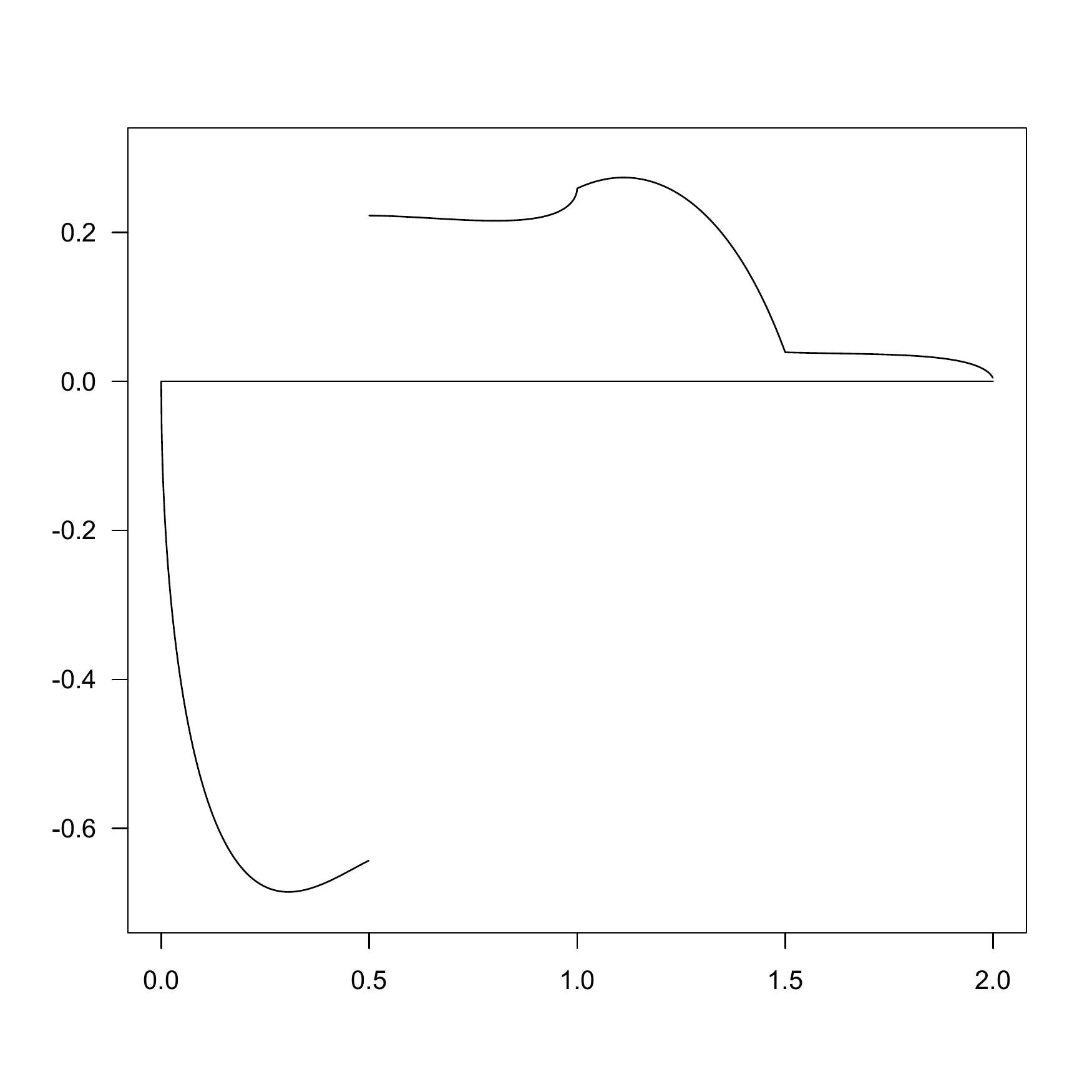}
\end{center}
\caption{The function $\th_{t,F_0}$, defined by  (\ref{score_decon_local}), for $t=0.5$, where $F_0$ is the uniform distribution function on $[0,1]$ and $g$ the elbow density $g(x)=2(1-x)1_{[0,1]}(x)$.}
\label{fig:theta_local}
\end{figure}

The corresponding efficient influence function $\th_{t,F_0}$, defined by
\begin{equation}
\label{score_decon_local}
\th_{t,F_0}(z)=\frac{\int_{u=0}^z g(z-u)\,d\f_{t,F_0}(u)}{h_0(z)}\,,
\end{equation}
is given in Figure \ref{fig:theta_local}.

\begin{figure}[!ht]
\begin{center}
\includegraphics[scale=0.5]{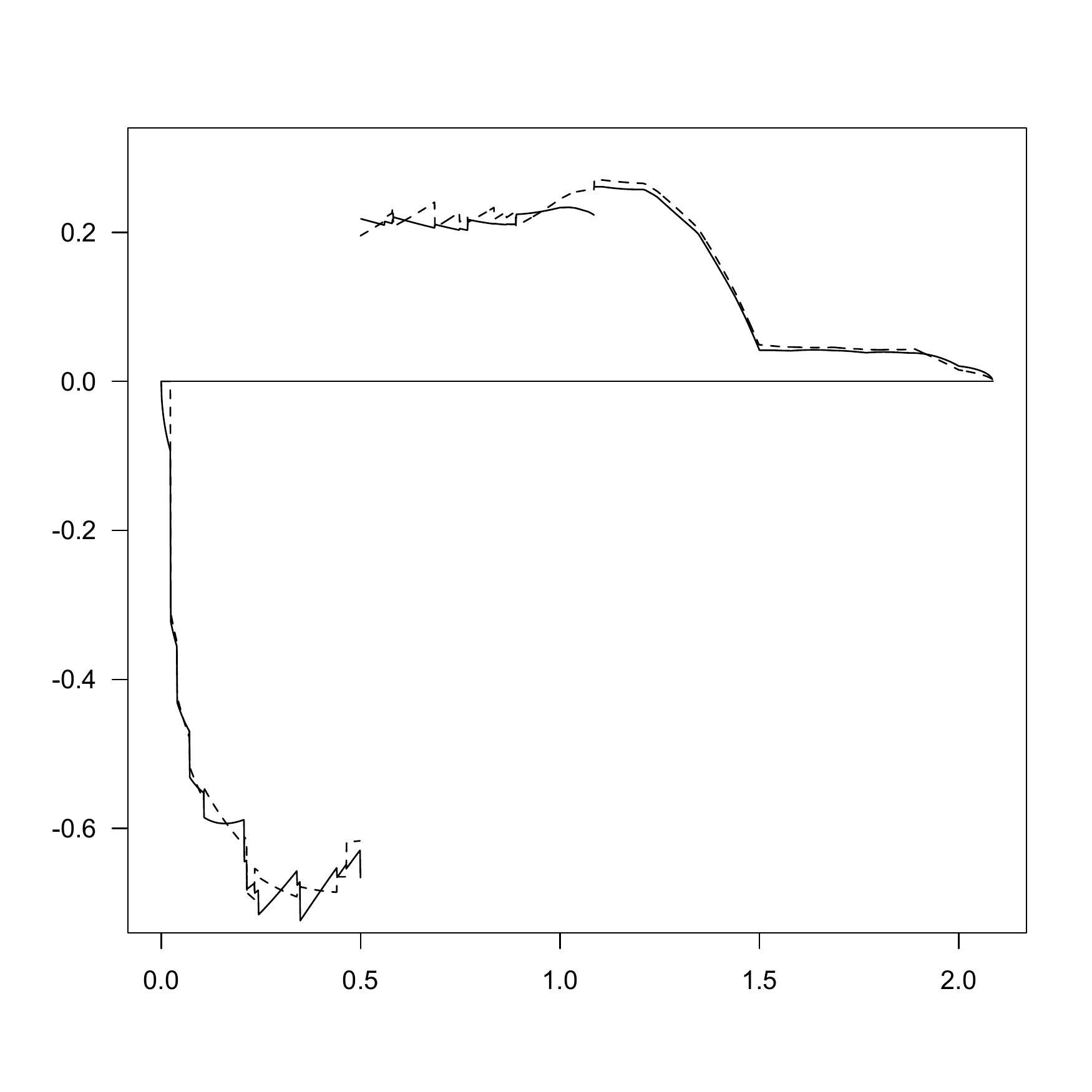}
\end{center}
\caption{The function $\th_{t,\hat F_n}\{\sqrt{\hat h_n}+\sqrt{h_0}\}/2$ (solid), defined by (\ref{thetaFn_representation_MLE_decon}),  and the function $\bar\th_{t,\hat F_n}\{\sqrt{\hat h_n}+\sqrt{h_0}\}/2$ (dashed), defined by (\ref{def_bar_theta}), for the same model as used in Figure \ref{fig:theta_local}, and for the same sample as used in Figure \ref{fig:phi_n_decon}; $\t_m\approx1.08617$.}
\label{fig:theta_n_decon_local}
\end{figure}

We can now proceed in a similar way as in the preceding section and define functions $\th_{t,\hat F_n}$ and $\bar\th_{t,\hat F_n}$ as in (\ref{thetaFn_representation_MLE_decon}) and (\ref{def_bar_theta}). The result is shown in Figure \ref{fig:theta_n_decon_local}. Note the analogy with what happened in section \ref{section:local_limit_IC}, where the proof of the convergence to the asymptotic distribution was also based on showing that a remainder term was of order $O_p(n^{-1/2})$ and asymptotically normal, by showing that an integral equation had a solution $\bar\f_{t,\hat F_n}$ which had a jump discontinuity at a point $t$ where we wanted to determine the local limit of the MLE, see (\ref{st-flour_inteq}). The difference with the present situation is that in that case a proof is available (in \cite{piet:96}), whereas for the deconvolution case the proof is still not completed.

It should be mentioned that the integral equation, determining $\f_{F}$ and $\th_{F}$, can be explicitly solved for exponential deconvolution. Note that we leave the compact support case here. Defining
$$
K_t(F)=\int_{x\in[0,t)}\{g(t-x)-g(0)\}\,dF(x),
$$
we get:
$$
\f_{t,\hat F_n}(x)=\left\{\begin{array}{lll}
-\{1+K_t(\hat F_n)\}\hat F_n(x),\,&x\in[\t_1,t),\\
-K_t(\hat F_n)\{\hat F_n(x)-1\},\,&x\in[t,\t_m],
\end{array}
\right.
$$
where $\t_1$ and $\t_m$ are the first and last point of jump of $\hat F_n$, respectively. Here we have the unusual situation that $\f_{t,\hat F_n}$ is almost absolutely continuous w.r.t. $\hat F_n$, which is only spoilt by the jump of $\f_{t,\hat F_n}$ at $t$. Moreover, we get:
$$
\th_{t,\hat F_n}(x)=-1_{[0,t)}(z)-K_t(\hat F_n).
$$
We have the following result.

\begin{theorem}
\label{th:expdec_smooth}
{\mbox{\rm($K_t(\hat F_n)$ is a smooth functional in exponential deconvolution.)}}
Let $X_1,\dots,X_n$ be a sample from a
continuous distribution function $F_0$, concentrated on $[0,\infty)$. Furthermore, let $Z_1,\dots,Z_n$ be a sample of observations of the type
$Z_i=X_i+Y_i$, where $Y_1,\dots,Y_n$ are independent of the $X_i$ and standard exponentially
distributed. Suppose that
$\hat F_n$ is the MLE of $F_0$ on the basis of the sample $Z_1,\dots,Z_n$. Then, for each point $t$ in the interior of the support of the distribution of the $X_i$:
\begin{equation}
\label{exp_limit}
\sqrt{n}\left\{K_t(\hat F_n)-K_t(F_0)\right\}=
\sqrt{n}\int_{x\in[0,t)}\bigl\{e^{-(t-x)}-1\bigr\}\,d\bigl(\hat F_n-F_0\bigr)(x)\stackrel{{\cal D}}\longrightarrow{\cal
N}(0,\s_t^2),
\end{equation}
where ${\cal N}(0,\s_t^2)$ denotes a normal distribution with first moment zero and variance $\s_t^2$, given by:
$$
\s_t^2=\int \th_{t,F_0}(z)^2\,dH_0(z),
$$
and where (the score function) $a_{t,F_0}$ is given by:
\begin{equation}
\label{a_exp0}
\th_{t,F_0}(z)=
-1_{[0,t)}(z)-K_t(F_0),\,z\ge0.
\end{equation}
\end{theorem}

Using this result, we get a new proof of the local limit result (\ref{local_limit_unif_expon}) above, since (\ref{smooth_local}) holds, and we get that the MLE is asymptotically equivalent to the ``toy estimator", obtained by doing one step of the iterative convex minorant algorithm, starting the iterations with the underlying distribution function $F_0$.

\section{Concluding remarks.} It is shown that further development of the local limit theory of the MLE for interval censoring and deconvolution crucially depends on getting grip on the associated integral equations. The same holds for the MSLE for interval censoring, which converges locally at a faster rate  than the MLE, under appropriate smoothing conditions. Probably similar results will follow for the MSLE for deconvolution, but this is not discussed above.

It is also shown that the MLE can be expected to be asymptotically efficient in the estimation of smooth functionals, a property it attains automatically, without any smoothing. For deconvolution this property will usually not hold automatically for the Fourier type estimators which are commonly applied in this situation, using kernel estimators with fixed bandwidths in the Fourier domain. The theory also produces answers to the often posed question: ``Why maximum likelihood?". One answer is: it produces automatically efficient estimates, in contrast with other methods. Nowadays, these estimates also can be easily computed, for example using the support reduction algorithm of \cite{piet_geurt_jon:08} or some hybrid form of the EM algorithm, combined with the iterative convex minorant algorithm, as was used for computing the MSLE in \cite{piet:12b}.

\section*{Acknowledgements}
I want to thank Anirban DasGupta for the invitation to write this article for the Journal of Statistical Planning and Inference.

\bibliographystyle{amsplain}
\bibliography{nonpar_int}
\end{document}